\documentclass{article}

\usepackage{arxiv}

\usepackage[utf8]{inputenc}
\usepackage[T1]{fontenc}
\usepackage{hyperref}
\usepackage{url}
\usepackage{doi}

\usepackage[english]{babel}
\usepackage{amsmath}
\usepackage{amsfonts}
\usepackage{graphicx}
\usepackage{natbib}
\usepackage{enumitem}

\usepackage{amssymb}

\usepackage{color}
\usepackage{multirow}
\usepackage{mathtools}
\usepackage{pdflscape}
\usepackage{tabularx} %for breaking text in tables 
\usepackage{mathrsfs} %for \matchscr font style
\usepackage{float} %For positioning the floats as H
\usepackage{booktabs} 
\usepackage{setspace}
\usepackage{threeparttable}
\usepackage{longtable}
\usepackage{subfigure}% Support for small, `sub' figures and tables

\usepackage[algoruled,vlined,linesnumbered,norelsize]{algorithm2e}
\SetKwComment{Comment}{\:\:\:\:$\triangleright$\ }{}
\allowdisplaybreaks

\newcommand{\cJ}{\mathcal{J}}
\newcommand{\cO}{\mathcal{O}}
\newcommand{\cM}{\mathcal{M}}

\newcommand{\cF}{\mathcal{F}}

\newcommand{\cU}{\mathcal{U}}
\newcommand{\cA}{\mathcal{A}}

\newcommand{\NP}{$\mathscr{N\mspace{-8mu}P}$}

\newcommand{\ihat}{{\hat{\imath}}}

\allowdisplaybreaks

\title{Iterated greedy algorithms for a complex parallel machine scheduling problem}

\author{
\href{https://orcid.org/0000-0003-0072-9793}{\includegraphics[scale=0.06]{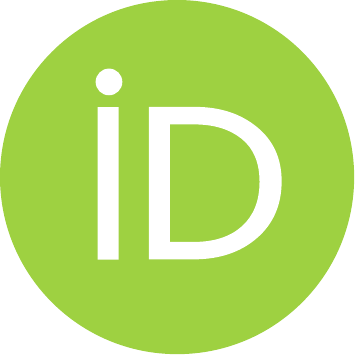}\hspace{1mm}Davi Mecler} \\
Departamento de Engenharia Industrial \\
Pontif\'{\i}cia Universidade Cat\'olica do Rio de Janeiro \\
\texttt{davizm@tecgraf.puc-rio.br} \\

\And
\href{https://orcid.org/0000-0003-3042-9249}{\includegraphics[scale=0.06]{orcid.pdf}\hspace{1mm}Victor Abu-Marrul} \\
Departamento de Engenharia Industrial \\
Pontif\'{\i}cia Universidade Cat\'olica do Rio de Janeiro \\
\texttt{victorabu@aluno.puc-rio.br} \\
\And
\href{https://orcid.org/0000-0001-5715-149X}{\includegraphics[scale=0.06]{orcid.pdf}\hspace{1mm}Rafael Martinelli} \\
Departamento de Engenharia Industrial \\
Pontif\'{\i}cia Universidade Cat\'olica do Rio de Janeiro \\
\texttt{martinelli@puc-rio.br} \\
\And
\phantom{xxxxxxxxxxxxxx}\href{https://orcid.org/0000-0002-8862-7603}{\includegraphics[scale=0.06]{orcid.pdf}\hspace{1mm}Arild Hoff} \phantom{xxxxxxxxxxxxxx}\\
Faculty of Logistics, \\
Molde University College \\
\texttt{arild.hoff@himolde.no}
}

\date{}

\renewcommand{\headeright}{A Preprint}
\renewcommand{\undertitle}{A Preprint}
\renewcommand{\shorttitle}{Iterated greedy algorithms for a complex parallel machine scheduling problem}

\hypersetup{
pdftitle={Iterated greedy algorithms for a complex parallel machine scheduling problem},
pdfsubject={math.OC},
pdfauthor={Davi Mecler, Victor Abu-Marrul, Rafael Martinelli, Arild Hoff},
pdfkeywords={Iterated greedy algorithms for a complex parallel machine scheduling problem},
}

\begin{document}
\maketitle

\begin{abstract}
This paper addresses a complex parallel machine scheduling problem with jobs divided into operations and operations grouped in families. Non-anticipatory family setup times are held at the beginning of each batch, defined by the combination of one setup-time and a sequence of operations from a unique family. Other aspects are also considered in the problem, such as release dates for operations and machines, operation's sizes, and machine's eligibility and capacity. We consider item availability to define the completion times of the operations within the batches, to minimize the total weighted completion time of jobs. We developed Iterated Greedy~(IG) algorithms combining destroy and repair operators with a Random Variable Neighborhood Descent~(RVND) local search procedure, using four neighborhood structures to solve the problem. The best algorithm variant outperforms the current literature methods for the problem, in terms of average deviation for the best solutions and computational times, in a known benchmark set of 72 instances. New upper bounds are also provided for some instances within this set. Besides, computational experiments are conducted to evaluate the proposed methods' performance in a more challenging set of instances introduced in this work. Two IG variants using a greedy repair operator showed superior performance with more than 70\% of the best solutions found uniquely by these variants. Despite the simplicity, the method using the most common destruction and repair operators presented the best results in different evaluated criteria, highlighting its potential and applicability in solving a complex machine scheduling problem.
\end{abstract}

\keywords{Metaheuristics \and Machine Scheduling \and Family Scheduling \and Batching \and Setup Times}

\section{Introduction}
\label{sec:introduction}

Parallel machine scheduling problems have been extensively studied since the seminal work of~\cite{mcnaughton1959scheduling}, with many applications in theoretical and real-life problems~\citep{wang2020identical,wang-2009}. The application that motivated this work is a real-life ship scheduling problem in the offshore oil and gas industry, with the machines representing an available fleet of vessels~\citep{AbuMArrul2020, abu2020matheuristics,cunhails}. The problem considers several aspects, such as non-anticipatory family setup times, eligibility constraints, split jobs, machine capacity, and others, making it more complex and challenging. Its characteristics relate to other parallel machine scheduling problems in the literature, such as serial batch scheduling problems, family scheduling problems, concurrent open shop scheduling problems, order scheduling problems, and some variants of the classical identical parallel machine scheduling problem.

More formally, there are a set of $n$ jobs, $\cal{J}$~=~\{$J_1$,...,~$J_n$\}, a set of $o$ operations, $\cal{O}$~=~\{$O_1$,...,~$O_{o}$\}, and a set of $m$ identical machines in parallel, $\cal{M}$~=~\{$M_1$,...,~$M_m$\}. Each job $J_j \in \cJ$ is composed by a subset $\cO_j \subseteq \cal{O}$ of operations, which can be processed by different machines simultaneously. All component operations must be finished to complete a job. In this problem, differently from job-shops, flow-shops, or open-shops, one operation might be associated with more than one job. Thus, we use $\cJ_i \subseteq \cJ$ to identify the subset of jobs that an operation $O_i$ is associated with. Operations are grouped by similarity, with each one belonging to exactly one family $g \in \cF$~=~\{$F_1$,...,~$F_{f}$\}. For each operation $O_i \in \cO$, its processing time $p_i$, release date $r_i$, and size $l_i$ are given. Moreover, a subset $\cM_i \subseteq \cM$ of eligible machines to execute operation $O_i$, without preemption, is known in advance. For each machine $M_k \in \cM$ is given a capacity $q_k$ and a release date $r_k$. When assigning operations to machines, a family setup time $s_g$ must be incurred whenever a machine changes the execution of operations from different families and before the first operation on each machine. The composition of a family setup time with a sequence of operations is called \textit{Batch}. The total size of all operations in a batch cannot exceed the assigned machine capacity. Finally, setup times are non-anticipatory since the starting time of a batch is restricted by the release date of the operations that composes it. The completion time $C_j$ of a job~$J_j$ is given by the maximum completion time among its component operations, that is, $C_j = \max\{C_i:O_i \in \cO_j\}$, where $C_i$ is the completion time of operation~$O_i$. The objective is to minimize the Total Weighted Completion Time~(TWCT) of jobs, given by $\sum_{j \in \cJ} w_jC_j$, where $w_j$ is the weight of job~$J_j$. 

In this problem, decision-makers need to not only decide the \textit{allocation} and \textit{sequencing} of operations in the machines but also how to \textit{compose batches}. The last part is critical since a sequence of operations of the same family in a machine can form batches in different ways due to capacity constraints. Furthermore, in some cases, splitting a sequence of operations of the same family into more than one batch may be a better choice, even if it does not violate the machine's capacity. This assumption considers the non-anticipatory setup aspect that may generate idleness in a given schedule, delaying processing some operations. 

\begin{figure}[htbp!] 
	\centerline{\includegraphics[scale=0.215, trim = 10 0 0 0, clip]{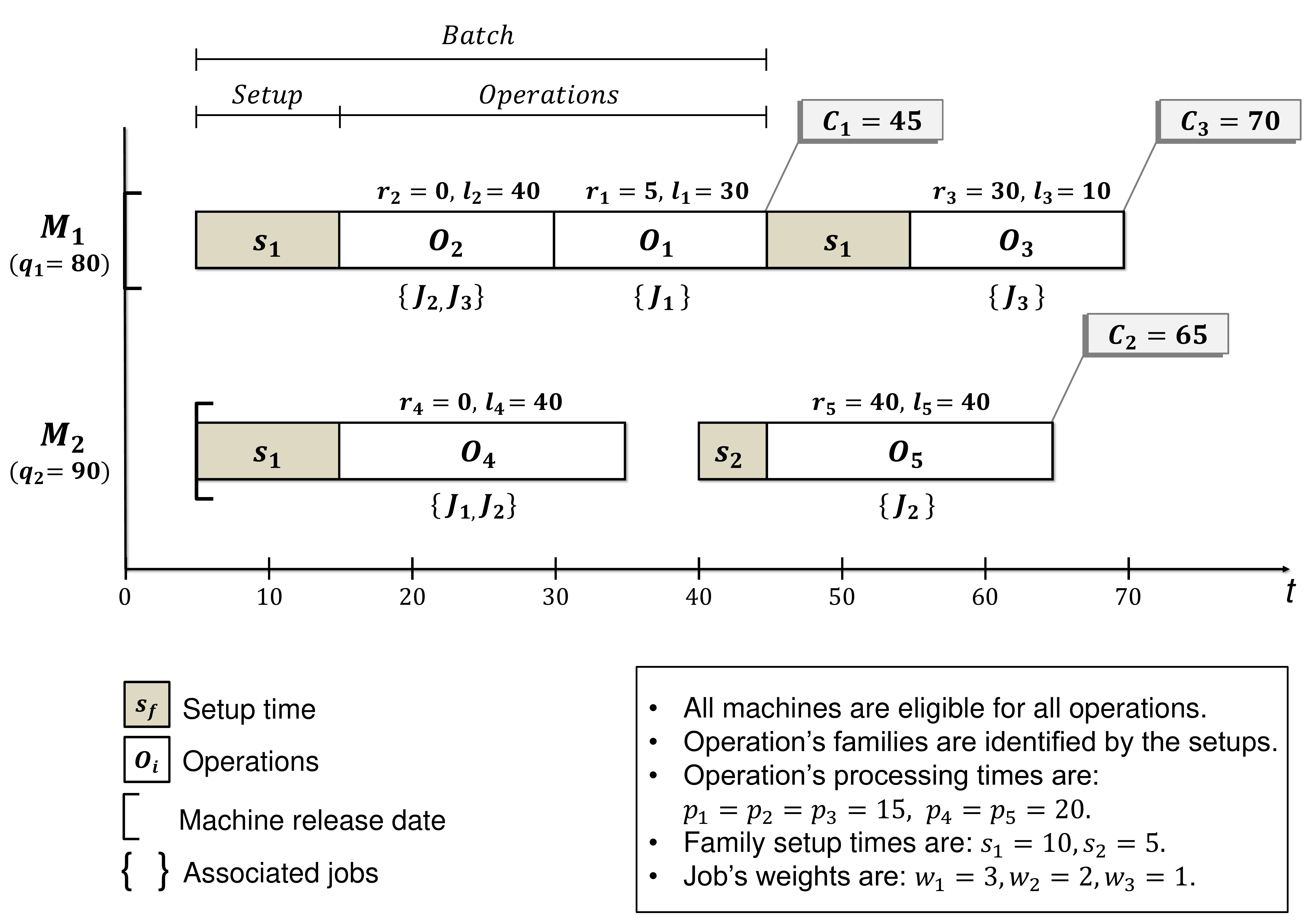}}
	\caption{Scheduling example with 5 operations, 3 Jobs and 2 machines. \label{fig:scheduling-example-one}}
\end{figure}

In Figure~\ref{fig:scheduling-example-one}, a small example, with 5 operations, 3 jobs, and 2 machines, is depicted to help the readers understanding the problem. At the top of each scheduled operation, we show its release date and size, and, below, we indicate its associated jobs. Setup times define the beginning of a batch, as shown in the top of the first batch on machine $M_1$. Each machine's capacity is shown below its name, and other characteristics are displayed inside the box below the schedule. Moreover, we highlight the completion times of the jobs in the schedule.

In the example given, two batches are created in each machine, respecting their capacities. Note that a single batch in machine $M_1$ would lead to a higher cost solution, since the unique batch would have to start its processing at $ t = 30 $ (highest release date between operations $O_1$, $O_2$, and $O_3$), respecting the non-anticipatory consideration. One can note that the non-anticipatory consideration generated idleness in both machines and that the first batch on machine $M_2$ respects the machine's release date. The Total Weighted Completion Time~(TWCT) of the example schedule, following the jobs index order $(J_1,J_2,J_{3})$, is given by: TWCT = (3 $\times$ 45) + (2 $\times$ 65) + (1 $\times$ 70) = 335.

In this paper, we extend the works of \citet{AbuMArrul2020, abu2020matheuristics} by proposing four variants of Iterated Greedy~(IG) algorithms, combined with a Randomized Variable Neighborhood Descent~(RVND) local search, to solve the described problem, testing it in two sets of benchmark instances. The paper's main contributions are: (1)~The development of IG approaches to solve a complex parallel machine scheduling problem, useful for similar machine scheduling problems; (2)~Improvement of the state of the art algorithm for the addressed problem; (3)~The development of a new set of instances, combining several concepts from distinct machine scheduling problems, making it more interesting to address; (4)~The combination of diversification and intensification strategies to improve the algorithms' performance; (5)~Solving a problem with a realistic ship scheduling background, highlighting the applicability of IG algorithms to real-life problems.

The outline of this paper is structured as follows. The problem's motivation regarding the machine scheduling theory is discussed, and some related works are presented in Section~\ref{sec:motivation}. Section~\ref{sec:ig_approach} details the developed algorithms' structure. Computational experiments are presented and discussed in Section~\ref{sec:experiments}. Finally, Section~\ref{sec:conclusion} concludes the paper, also giving some perspectives for future research regarding the problem.

\section{Motivation and Related Works}
\label{sec:motivation}

As stated in Section~\ref{sec:introduction}, the paper's initial motivation comes from a ship scheduling problem related to the offshore oil and gas industry. In the problem, a set of vessels, analogous to the machines, must perform operations in sub-sea oil wells to anticipate processing the most productive ones. Regarding this problem, referred to in the literature as the \textit{Pipe Laying Support Vessel Scheduling Problem}~(PLSVSP), \citet{AbuMArrul2020} developed three mathematical formulations, testing it in a set of instances, also generated by the authors, based on real data provided by a Brazilian company. \citet{abu2020matheuristics} developed six matheuristics, combining a constructive heuristic with two batch scheduling formulations. The authors provided new best solutions for several instances in the benchmark of \citet{AbuMArrul2020}, reducing the computational time compared with the pure mathematical formulations. 

From the best of our knowledge, no other work in the scheduling literature deals simultaneously with all the studied problem aspects. However, some works have certain similarities with the problem. \citet{Gerodimos-1999} approached a single machine scheduling problem with jobs divided into operations from different families, including sequence-independent setup times. The authors provided a complexity analysis of the problem within three objectives to minimize: maximum lateness, weighted number of late jobs, and total completion time. \citet{mosheiov-2008} proposed a linear time algorithm to minimize makespan and flow-time in a problem in which jobs are split into operations with identical processing times. Non-anticipatory setup times are incurred between batches, with their sizes limited by the machine's capacities. In \citet{wang-2009}, jobs are divided into orders that are organized into families. Sequence-dependent family setup times are considered, and jobs are delivered in batches with limited capacity. The authors proposed heuristic algorithms to minimize the weighted sum of the last arrival time of jobs. \citet{tahar-2009} proposed a linear programming approach to minimize the makespan in a parallel machine scheduling problem with job sections and sequence-dependent setup times between sections from different jobs. A sequence-dependent family setup time marks the starting process of a batch. In \citet{shen-2012}, a sequence-dependent family setup time marks the starting process of a batch. They developed a tabu search algorithm to minimize the makespan in a problem with jobs organized into families and composed by operations. In \citet{kim-2018}, jobs are processed in parts, organized in families, on a set of parallel machines to minimize the makespan. Setup times are incurred between processing parts from different families. \citet{Rocholl-2019} considered a problem with jobs divided into orders from the same family, with equal processing times. The authors proposed a mixed-integer linear problem model to minimize both earliness and tardiness, respecting the machine's capacities. \citet{lin2011} addressed a concurrent open shop problem, where jobs are partitioned into operations that are eligible to only a specific machine. The authors provided polynomial-time algorithms to minimize the maximum lateness, the weighted number of tardy jobs, and the total weighted completion time.

Moreover, due to the addressed problem complexity, it can be seen as a generalization of other well-known machine scheduling problems. This relation highlights the relevance of the problem regarding the scheduling literature. Still, it indicates that the algorithms presented in this work may be used to solve all simplified variants of the problem. In the following sections, we list four machine scheduling problems that can be seen as simplifications of the addressed problem.

\subsection{Identical Parallel Machine Scheduling Problem}
\label{sec:identical_p_machine}

The first related problem is the classical \textit{Identical Parallel Machine Scheduling Problem}~(IPMSP) that can be achieved by simplifying four aspects: (1) Every operation is associated to a single job; (2) Every job is composed by a single operation; (3) Every operation belongs to the same family; (4) Setup time duration is set as zero for all families. Note that, with these modifications, capacity constraints are disregarded, since no setup times are considered. Furthermore, each operation completion time will represent the completion time of its unique associated job. Observe that release dates can be considered or not, since it limits the starting of a batch with a single operation and zero duration. Thus, it will limit the operation itself. Formally, there is a set of $n$ jobs, $\cal{J}$~=~\{$J_1$,...,~$J_n$\}, to be processed by a set of $m$ identical machines in parallel, $\cal{M}$~=~\{$M_1$,...,~$M_m$\}. Each job~$J_j$ has a processing time $p_j$, might have a release date~$r_j$, and might be limited to be processed by a subset~ $\cM_j \subseteq \cM$ of eligible machines. The Identical Parallel Machine Scheduling Problem has been proven to be \NP-hard \citep{mokotoff-2004}. Thus, we can assume that the addressed problem is \NP-hard according to the described simplification.

\citet{nessah2008} and \citet{ahmed2013} addressed the IPMSP with release dates. In the former, a branch-and-bound algorithm is developed to minimize the total weighted completion time. In the latter, a modified forward heuristic to minimize makespan is presented. \citet{xu2014} proposed a column generation approach to minimize both makespan and total weighted completion time. \citet{dellamico2008}, \citet{alharkan2018}, and \citet{dellacroce2018} developed heuristic algorithms to deal with the IPMSP to minimize the makespan, while \citet{dellamico2005} solution approaches are based on exact methods. \citet{haouari2006} developed both heuristics and exact methods considering the same objective of minimizing the makespan.

\subsection{Family Scheduling Problem}
\label{sec:family_p_machine}

The second related problem is the \textit{Family Scheduling Problem}~(FSP), achieved with the following simplifications: (1)~Every operation is associated with a unique job; (2)~Every job is composed by a single operation; (3)~Release dates are set as zero; (4)~Size of operations are set as zero.  Note that, with these modifications, family setup times are only considered when changing the execution from operations of different families. Furthermore, the completion time of one operation will represent the completion time of its single related job. If release dates are considered, setup times would be non-anticipatory and additional setups may be considered during the optimization to avoid delays in starting processing some operations. Formally, there is a set $n$ of jobs, $\cal{J}$~=~\{$J_1$,...,~$J_n$\}, to be processed by a set of $m$ identical machines in parallel, $\cal{M}$~=~\{$M_1$,...,~$M_m$\}. Each job $J_j \in \cJ$ belongs to exactly one family $g \in \cF$~=~\{$\cF_1$,...,~$\cF_{f}$\}. Moreover, each job $J_j$ has a processing time $p_j$, might have a release date $r_j$, and might be limited to be processed by a subset $\cM_j \subseteq \cM$ of eligible machines. A family setup time $s_g$ is incurred before the first job on each machine or when a machine changes the execution of jobs from different families.

\citet{mehdizadeh2015}, \citet{liao2012}, \citet{bettayeb2008} and \citet{dunstall2005} addressed the FSP with parallel machines. The first work proposed a vibration-damping metaheuristic. The second and the third ones developed constructive heuristics, and the last introduced a branch-and-bound algorithm to minimize the total weighted completion time of jobs. \citet{lin2017} and \citet{Nazif2009} addressed the FSP on a single machine. The former considered two variants of the problem, with and without due dates, introducing a mixed integer programming model to minimize the total weighted completion time and maximum lateness. The latter work developed a genetic algorithm to minimize the total weighted completion time.

\subsection{Concurrent Open Shop Scheduling Problem}
\label{sec:open_shop}

The \textit{Concurrent Open Shop Scheduling Problem}~(COSSP), also referred to as the \textit{Order Scheduling Problem}~(OSP), is the third related problem. Three simplifications are needed to achieve it: (1) Every operation is associated with a single job; (2) Setup times are set as zero; (3) Each operation has a unique eligible machine to execute it. In this problem, a set of operations must be executed, each one in its pre-defined machine, to complete a job. Note that the variants with and without release dates can be solved since no setup times are considered in this problem. Formally, there are a set of $n$ independent jobs, $\cal{J}$~=~\{$J_1$,...,~$J_n$\}, and a set of $m$ identical machines in parallel, $\cal{M}$~=~\{$M_1$,...,~$M_m$\}. Jobs are composed of operations defined in the set $\cal{O}$~=~\{$O_1$,...,~$O_o$\}. The subset of operations that composes a job $J_j \in \cJ$ is given by $\cO_j \subseteq \cal{O}$, in such a way that $\bigcup_{j=1}^{n} \cO_{j} = \cO$. Each operation $O_i$ is eligible to be processed by a unique machine, it has a processing time $p_i$, and might have a release date $r_i$. All component operations must be executed to complete a job. 

\citet{Roemer2006} addressed the COSSP, providing a complexity analysis aiming to minimize the makespan. \citet{mastrolilli2010} proposed a combinatorial approximation algorithm to minimize the total weighted completion time, while \citet{Cheng2011} introduced a polynomial-time algorithm with the same objective. \citet{khuller2019} considered release dates, introducing an online scheduling framework, and minimizing the total weighted completion time.

\subsection{Serial Batch Scheduling Problem with Job Availability}
\label{sec:serial_batch}

The last analogy regards to a \textit{Serial Batch Scheduling problem with Job availability}, and it can be achieved with four simplifications: (1) Every operation is associated with one job; (2) Every job is composed by a single operation; (3) Release dates are set as zero; (4) All operations belong to the same family. In this problem, a setup time with a fixed duration is incurred at the beginning of each batch, and a batch must respect the machine capacity. If we consider release dates, the setup times become non-anticipatory. Formally, there are a set of $n$ independent jobs, $\cal{J}$~=~\{$J_1$,...,~$J_n$\}, to be processed by a set of $m$ identical machines in parallel, $\cal{M}$~=~\{$M_1$,...,~$M_m$\}. Each job $J_j$ has a processing time $p_j$ and a size $l_j$. A fixed setup time $s$ must be incurred before the first job on each machine or when the machine capacity $q_k$ is reached. 

Serial batch scheduling problems with job availability are rarely found in the scheduling literature. \citet{gerodimos2001scheduling} addressed the problem with a single-machine and job's due dates. They introduced a pseudo-polynomial time algorithm to minimize three objective functions: total weighted completion time, maximum lateness, and the number of tardy jobs.

\section{Iterated Greedy Algorithm}
\label{sec:ig_approach}

In this section, we describe the structure of the proposed IG algorithms. The IG approach~\citep{ruiz2007simple} combines an intensification step, given by a local search procedure to achieve local optimal solutions, with destroying and repair phases to diversify solutions and avoid the method being stuck. Our IG algorithm uses a Random Variable Neighborhood Descent~(RVND) in the local search phase. The RVND strategy \citep{subramanian2017} picks neighborhoods up randomly from a pool instead of running them in a deterministic pre-defined sequence. We also included a simulated annealing criterion, a strategy to accept infeasible solutions to add more diversification to the method during its execution, and a solution restore step to intensify the search around good feasible solutions. 

\begin{figure}[h!]
\centering
    \includegraphics[scale=0.2, trim = 0 0 0 0,clip]{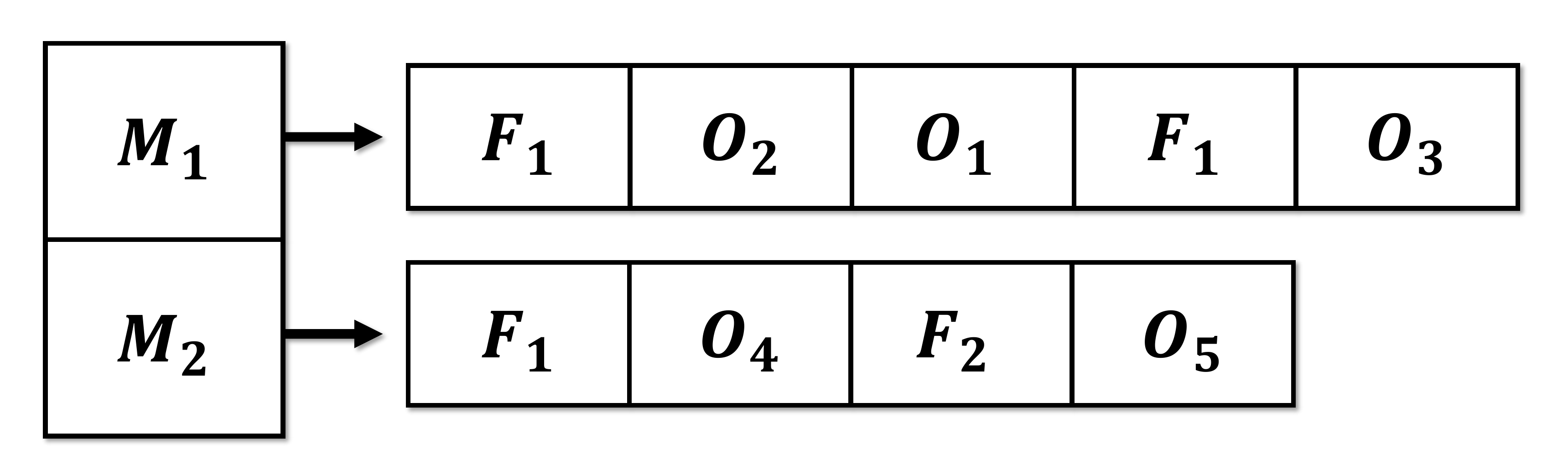}
    \caption{Solution representation regarding the schedule example depicted in Figure~\ref{fig:scheduling-example-one}}
    \label{fig:sol_rep}
\end{figure}

The general pseudo-code of the IG approach is shown in Algorithm~\ref{alg:ils}. The method considers $\mathtt{s}=(\mathtt{s}_1,\ldots,\mathtt{s}_m)$ as a solution composed by a list of schedules for the $m$~machines. Each schedule is a permutation of families (indicating a setup time) and operations. $\mathtt{s}$ indicates the current solution on each iteration and $\mathtt{s^*}$ the best solution reached so far in the execution. And, function $f(.)$ returns the total weighted completion time of an input solution. In Figure~\ref{fig:sol_rep}, the solution representation is presented, considering the schedule example depicted in Figure~\ref{fig:scheduling-example-one}.

\begin{algorithm}[!htb] %[htb!]
\footnotesize
	%\setstretch{1.2}
	\SetAlgoLined
	\SetKwInOut{Input}{Input}
  \SetKwInOut{Output}{Output}

\Input{Number of iterations ($\eta$); Restore parameter ($\lambda$).}
  \Output{The best solution found $\mathtt{s^*}$ as a list of schedules for each machine.}
  \BlankLine
    	
    	$\Omega \gets \lceil \lambda  \eta \rceil$;
    	
        $\mathtt{s} \gets \texttt{Constructive()}$;\Comment*[f]{Construct the initial solution}
        
    	$\mathtt{s} \gets$ \texttt{RVND($\mathtt{s}$)};\Comment*[f]{Local Search}
    	
    	$\mathtt{s^*} \gets \mathtt{s}$;\Comment*[f]{Initialize the best overall solution}
    	
    	$\omega \gets 0$; 
    	
	   	\For{$\eta$ iterations}{

	   	$\mathtt{s'} \gets$ \texttt{Destroy($\mathtt{s}$)}; 
	   	
	   	$\mathtt{s'} \gets$ \texttt{Repair($\mathtt{s'}$)};
		
	    $\mathtt{s'} \gets$ \texttt{RVND($\mathtt{s'}$)};
	    
	    $\omega \gets \omega + 1$;
	    
	    $is\_feasible \gets$ \texttt{Feasible($\mathtt{s'}$)};\Comment*[f]{Infeasibility strategy}
	    
	   \If(\Comment*[f]{Acceptance Criterion}){$f(\mathtt{s'}) < f(\mathtt{s})$ {\normalfont \textbf{or} \texttt{Accept($\mathtt{s'},\mathtt{s}$)}}}{
	    
	    $\mathtt{s} \gets \mathtt{s'}$;

		\If(\Comment*[f]{Improvement check}){$f(\mathtt{s}) < f(\mathtt{s^*})$ {\normalfont \textbf{and}}  $is\_feasible$}{

		$\mathtt{s^{*}} \gets \mathtt{s}$;
		
		$\omega \gets 0$;

		}
		}
		
		          \If(\Comment*[f]{Solution restore}){$\omega$ = $\Omega$}{
		
 		$\mathtt{s} \gets \mathtt{s^{*}}$;
		
 		$\omega \gets 0$;
		
 		}

	    }
	
	\textbf{return} $\mathtt{s^*}$;
		
	\caption{Iterated Greedy Algorithm 
	} 
	\label{alg:ils}
\end{algorithm}

Algorithm~\ref{alg:ils} starts by computing the number of consecutive iterations without improvement $\Omega$ to restore the current solution~(Line~1). Then, a solution $\mathtt{s}$ is built by employing the \texttt{WMCT-WAVGA} constructive heuristic~(Section~\ref{sec:constructive})\citep{abu2020matheuristics}~(Line~2). After running the local search~(Line~3), the best overall solution is initialized~(Line~4), and also the counter $\omega$ of consecutive iterations without improvement~(Line~5). The main loop~(Lines~6--23) runs for $\eta$ iterations and consists of destroying and repairing the current solution $\mathtt{s}$, generating solution $\mathtt{s'}$~(Lines~7--8), and applying a local search in $\mathtt{s'}$ to achieve a local optima solution~(Line~9). Then, the feasibility of the local optima solution is verified~(Line~11). The algorithm updates the current solution $\mathtt{s}$ if the cost of $\mathtt{s}'$ is lower than the cost of $\mathtt{s}$ or if the acceptance criterion is met~(Lines~12--13). If the updated current solution $\mathtt{s}$ is feasible and its cost improves the cost of solution $\mathtt{s^*}$, $\mathtt{s^*}$ is replaced by $\mathtt{s}$ and the counter $\omega$ of iterations without improvement is reinitialized~(Lines~14--17). The algorithm replaces $\mathtt{s}$ by $\mathtt{s^*}$ after $\Omega$ consecutive iterations without improvement~(Lines~19--22). The method returns the best overall solution $\mathtt{s^*}$~(Line~24). 

Next sections detail the constructive heuristic, the RVND local search with the considered neighborhoods, the destroy and repair operators, the acceptance criterion, and the infeasibility strategy. 

\subsection{Constructive Heuristic}
\label{sec:constructive}

The \texttt{Constructive} procedure (Line 2 of Algorithm~\ref{alg:ils}) builds solutions iteratively by using the \texttt{WMCT-WAVGA} heuristic~\citep{abu2020matheuristics}, shown in Algorithm~\ref{alg.constructive}, assigning  operations sequentially to the machines. The heuristic combines an adaptive rule to estimate the operations' weights with a weighted minimum completion time dispatching rule to select the operations to schedule. The selected operation is scheduled on the eligible machine, which minimizes the partial schedule's total weighted completion time.

\renewcommand{\gets}{\leftarrow}
\begin{algorithm}[htb]
	    \renewcommand{\gets}{\leftarrow}
	%\setstretch{0.8}
	\footnotesize
	\SetAlgoLined
		\SetAlgoLined
	\SetKwInOut{Input}{Input}
  \SetKwInOut{Output}{Output}

\Input{Problem data.}
  \Output{A solution $\mathtt{s}$.}
  \BlankLine

			$C_k \gets r_k$, $S_k \gets r_k$, $L_k \gets 0$, $F_k \gets 0$, $\cA_k \gets \emptyset$, $\mathtt{s}_{k} \gets \emptyset$, \:\: $\forall M_k \in \cM$;

			$C_i \gets \infty, \:\: \forall O_i \in \cO$;

			$\mathcal{U} \gets \cO$;

		\While{$\cU \neq \varnothing$}{
		
		$w_i \gets \sum\limits_{j \in \mathcal{N}_i}\frac{w_j}{|\cU_j|}$, $T_i \gets \min\limits_{k \in \cM_i}C_k$, $\pi_i = \frac{w_i}{\max(T_i, r_i) + p_i + s_{g_i}}$, \:\:$\forall O_i \in \mathcal{U}$;
		
		Select operation $O_{i^*} \in \mathcal{U}$ that maximizes $\pi_i$;

		$\Delta_{i^* k} \gets \max(0, r_{i^*} - S_k),  \:\: \forall M_k \in \cM_{i^*}$;

		$C^{CB}_{i^* k} \gets C_k + \Delta_{i^* k} + p_{i^*}, \:\: \forall M_k \in \cM_{i^*}$;

		$C^{NB}_{i^* k} \gets \max(r_{i^*}, C_k)+s_{g_{i^*}} + p_{i^*},  \:\: \forall M_k \in \cM_{i^*}$;
		
		$ \mathcal{CB} \gets \bigg\{ cb_{i^* k} = w_{i^*} C^{CB}_{i^* k} + \sum\limits_{\ihat \in \mathcal{B}_k }{w_\ihat} \Delta_{i^* k} \:\: \big| \:\: k \in \cM_{i^*}, \: F_k=g_{i^*}, \: L_k+l_{i^*} \leq q_k \bigg\}$;
		    
	    $ \mathcal{NB} \gets \bigg\{nb_{i^* k} =w_{i^*} C^{NB}_{i^* k}\:\: \big| \:\: M_k \in \cM_{i^*}\bigg\}$;
	    
        $b_{min} \gets \min\{b :  b \in (\mathcal{CB} \cup \mathcal{NB})\}$;
        
		Select $M_{k^*}$ corresponding to $b_{min}$;
	
		\eIf{$b_{min} \in \mathcal{CB}$}{

		$S_{k^*} \gets \max(r_{i^*}, S_{k^*})$; $C_{k^*} \gets C_{i^* k^*}^{CB}$;
			
		$L_{k^*} \gets L_{k^*} + l_{i^*}$; $\cA_{k^*} \gets \cA_{k^*} \cup \{O_{i^*}\}$;

		}{    
		
		$S_{k^*} \gets \max(r_{i^*}, C_{k^*})$; $C_{k^*} \gets C_{i^* k^*}^{NB}$;
		
		$L_{k^*} \gets l_{i^*}$; $\cA_{k^*} \gets \{O_{i^*}\}$;
		
		$\mathtt{s}_{k^*} \gets  \mathtt{s}_{k^*} \cup \{g_{i^*}\}$;
		
		}
		
		$\mathtt{s}_{k^*} \gets  \mathtt{s}_{k^*} \cup \{O_{i^*}\}$; $F_{k^*} \gets g_{i^*}$; $C_{i^*} \gets C_{k^*}$; $\mathcal{U} \gets \mathcal{U} \setminus \{O_{i^*}\}$;
		}
		\textbf{return} $\mathtt{s}$;

	\caption{\texttt{WMCT-WAVGA \citep{abu2020matheuristics}}} 
	\label{alg.constructive}
\end{algorithm}

The algorithm starts by initializing $S_k$, $L_k$, $F_k$, and $\cA_k$, used to store the starting time, cumulative load, family, and set of assigned operations regarding the current batch on machine~$M_k$, respectively~(Line~1). Completion times of operations are initialized~(Line~2), and the set of unscheduled operations is created~(Line~3). The algorithm's main loop~(Lines~4--23) stops when all operations are assigned. Within the loop, operations are selected~(Line~6) according to a priority index computed in line~5. We use $g_i$ to indicate the family of operation~$O_i$. The delay in starting the current batch~($\Delta_{ik}$) on each machine~$M_k$, considering the insertion of the selected operation~$O_{i^{*}}$ is computed in line~7. The value of $\Delta_{ik}$ is used to compute the values of $C_{i^{*}k}^{CB}$ and $C_{i^{*}k}^{NB}$. They represent the completion time of the selected operation~$O_{i^*}$, if inserted in the current batch or in a new batch on each machine~$M_k$, respectively. Then, the algorithm identifies the sets $\mathcal{CB}$ and $\mathcal{NB}$ with the feasible assignments of the selected operation into the current batch or into a new batch, on each eligible machine~$M_k$, with the respective costs~(Lines~10--11). The best element among the sets is chosen~(Line~12), and the corresponding machine is identified~(Line~13). Finally, the selected operation is assigned to the selected machine, and the algorithm's variables and the schedule are updated~(Lines~14--22). The algorithm returns the final schedule $\mathtt{s}$~(Line~24).

\subsection{RVND Local Search}
\label{secsim:rvnd}

\begin{figure}[htb]
\centering
	\subfigure[Swap between operations of the same family.\label{fig:swap_same}]{%
		\includegraphics[width=1\textwidth, trim = 0 0 0 0, clip]{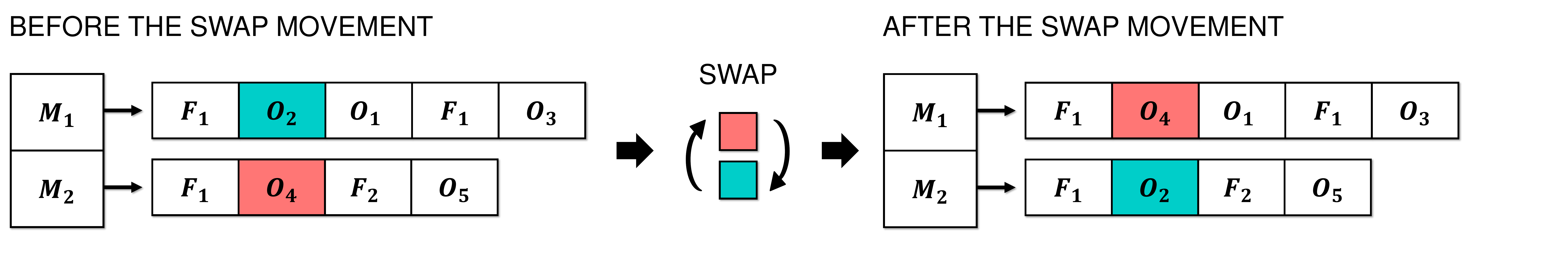}
	}
	\subfigure[Swap between operations of distinct families.\label{fig:swap_distinct}]{%
		\includegraphics[width=1\textwidth, trim = 0 0 0 0, clip]{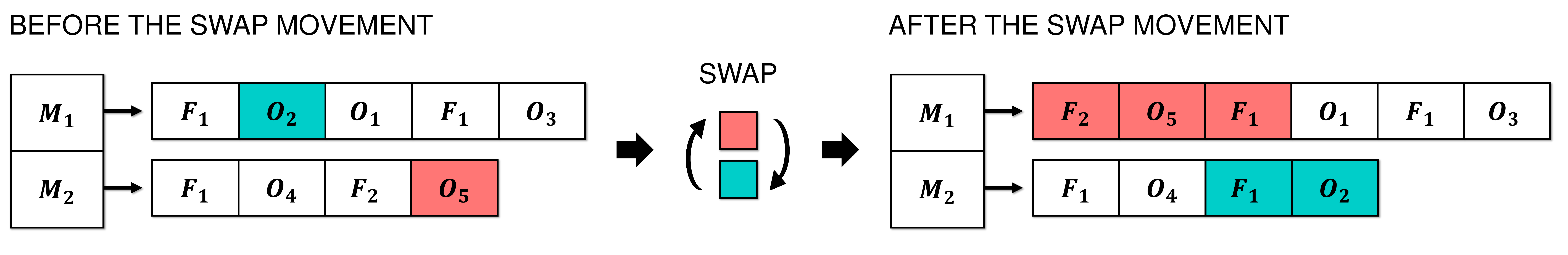}
	}
	\caption{Swap movement examples.}
	\label{fig:swap}
\end{figure}

\begin{figure}[htb]
\centering
	\subfigure[Relocate in a batch of the operation's family.\label{fig:relocate_same}]{%
		\includegraphics[width=1\textwidth, trim = 0 0 0 0, clip]{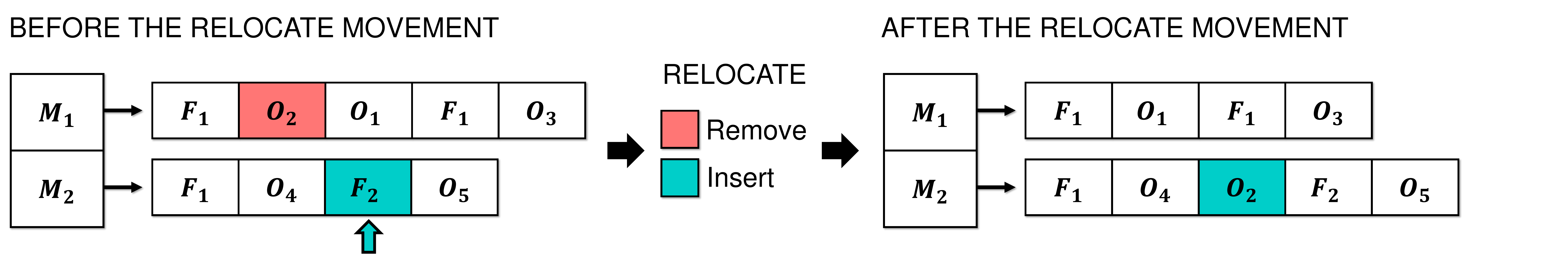}
	}
	\subfigure[Relocate in a batch of a distinct families.\label{fig:relocate_distinct}]{%
		\includegraphics[width=1\textwidth, trim = 0 0 0 0, clip]{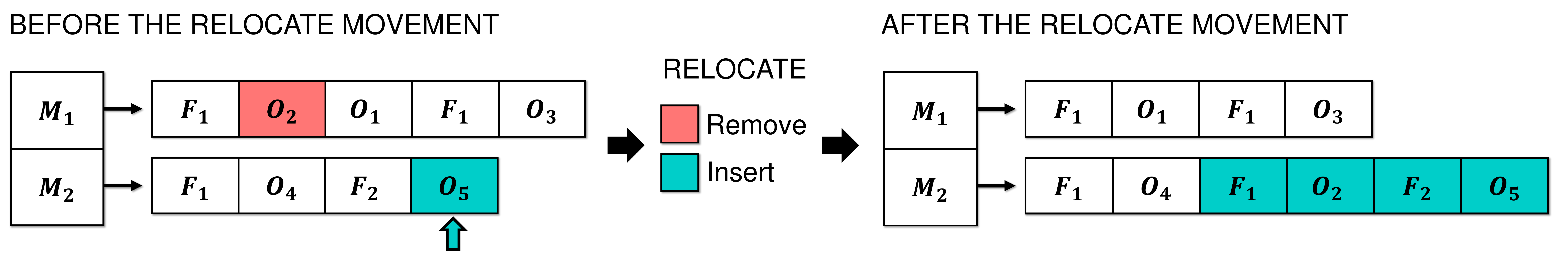}
	}
	\caption{Relocate movement examples.}
	\label{fig:relocate}
\end{figure}

\begin{figure}[htb]
\centering
	\subfigure[Split Batches movement.\label{fig:setup_insert}]{%
		\includegraphics[width=1\textwidth, trim = 0 0 0 0, clip]{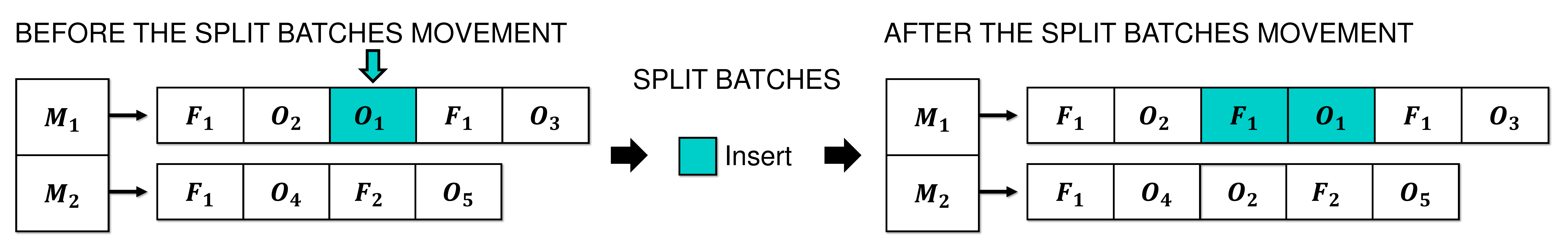}
	}
	\subfigure[Merge Batches movement.\label{fig:setup_removal}]{%
		\includegraphics[width=1\textwidth, trim = 0 0 0 0, clip]{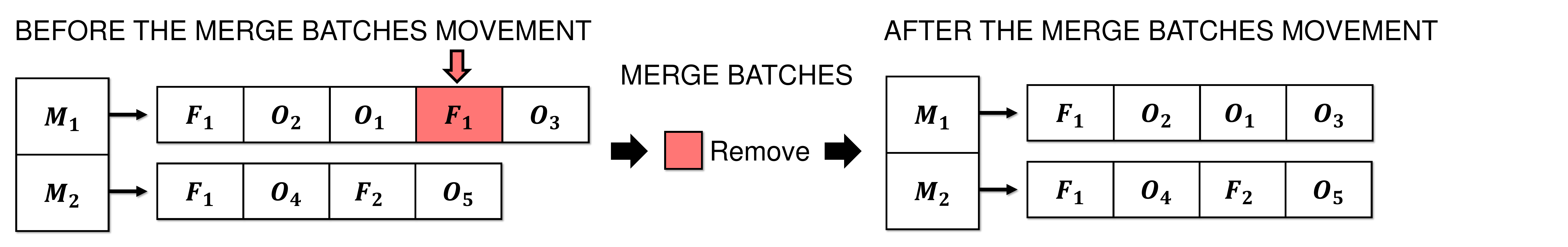}
	}
	\caption{Setup movement examples.}
	\label{fig:setup_movements}
\end{figure}

The \texttt{RVND} procedure~(Lines 3 and 9 of Algorithm~\ref{alg:ils}) takes a solution and returns a local optimum solution, considering four neighborhood structures, described as follows:

\begin{enumerate}

\item \texttt{Swap}: Select any two operations, assigned to the same or different machines, and exchange them on the schedule, considering the eligibility constraint. Figure~\ref{fig:swap} exemplifies two swap movements. The first between operations from the same family (Figure~\ref{fig:swap_same}), and the second between operations from distinct families (Figure~\ref{fig:swap_distinct}). Note that, when the movement generates batches with mixed families, new setup times are included and extra setup times are removed to make the solution feasible regarding the family constraint. 

\item \texttt{Relocate}: Select an operation assigned to any machine and insert it in a new position on the same or another machine, considering the eligibility constraint.  Figure~\ref{fig:relocate} exemplifies two relocate movements. The first inserts an operation in a batch of the same family (Figure~\ref{fig:relocate_same}), and the second inserts an operation in a batch of a distinct family (Figure~\ref{fig:relocate_distinct}). Again, one can note that setup times are included to ensure the feasibility of the solution regarding the family constraint. Also, if a movement creates a solution with two consecutive setup times, the extra setup time is removed.

\item \texttt{SplitBatches}: Select a batch and insert a setup time between two consecutive operations of the same family. The idea is to split a batch into two.
Figure~\ref{fig:setup_insert} depicts an example of a setup insertion movement. Note that, this movement is only possible in batches with more than one operation inside. 

\item \texttt{MergeBatches}: Select a setup time on any machine and remove it, if possible. A batch from the same family must precede the batch with the chosen setup time for the movement to be feasible. The idea is to combine batches from the same family scheduled in sequence. Figure~\ref{fig:setup_removal} depicts an example of a setup removal movement. Note that, this movement is only possible when the resulting solution is feasible regarding the family constraint.
\end{enumerate}

The \texttt{RVND}, shown in Algorithm~\ref{alg.rvnd}, starts by initializing the set~$\cal{L}$ with the indices of the considered neighborhoods in the local search procedure~(Line~1). Considering the proposed neighborhoods (1--\texttt{Swap}, 2--\texttt{Relocate}, 3--\texttt{SplitBatches}, 4--\texttt{MergeBatches}), the set of indices is $\cal{L}=\{$1, 2, 3, 4$\}$. The main loop is repeated until set $\cal{L}$ is empty~(Lines~2--17). The algorithm selects a neighborhood index $\ell$ at random~(Line 3) and tests each solution in the neighborhood $N_{\ell}$ until it finds the first improvement in the solution cost~(Lines~5--11). If an improved solution is found, the current solution is updated, and the algorithm resets the set $\cal{L}$ of neighborhood indices. If no improvement is found in a given neighborhood, its index is removed from set $\cal{L}$~(Lines~12--16). The algorithm returns the local optimal solution~(Line 18).

\begin{algorithm}[htb!]
    \renewcommand{\gets}{\leftarrow}
	%\setstretch{0.8}
	\footnotesize
	\SetAlgoLined
		\SetAlgoLined
	\SetKwInOut{Input}{Input}
  \SetKwInOut{Output}{Output}

\Input{A given solution $\mathtt{s}$.}
  \Output{A local optimal solution $\mathtt{s}$.}
  \BlankLine
	   
	     Initialize set $\cal{L}$ of neighborhood indices; 
	        
	        \While{$\cal{L} \neq \varnothing$} {
	        
	        $\ell \gets$\texttt{random}($\cal{L}$);
	        
	            $improved \gets$ \texttt{false};
	        
    	        \For{$\mathtt{s'} \in N_{\ell}(\mathtt{s})$} {
    	        
                      	            \If{$f(\mathtt{s'}) < f(\mathtt{s})$} {
    	            
    	                $\mathtt{s} \gets \mathtt{s'}$;
    	                
    	                $improved \gets$ \texttt{true};

    	                \texttt{break};
    	            
    	            }
    
    	      }
    	      \uIf{improved}{

    	      Initialize set $\cal{L}$ of neighborhood indices;
    	      }
    	      \Else{
    	      $\cal{L} \gets \cal{L} \setminus \{\ell\}$;
    	      }

	    }
	
	\textbf{return} $\mathtt{s}$;
		
	\caption{\texttt{RVND}} 
	\label{alg.rvnd}
\end{algorithm}

\subsection{Destroy and Repair Operators}
\label{sec:destroy_repair}

The \texttt{Destroy} procedure~(Line 7 of Algorithm~\ref{alg:ils}) removes $d = \lceil \varepsilon   o \rceil$ operations from a given solution $\mathtt{s}$, generating a partial solution with $o - d$ scheduled operations, where $\varepsilon \in [0, 1]$ is the destruction parameter. Extra setup times are also removed from the partial solution when necessary. The total weighted completion time of a partial solution only considers the scheduled operations, equivalent to defining all unscheduled operations' completion times to zero. It is worth mentioning that the batches' starting times are anticipated after removing the operations, whenever possible, to calculate the partial solution's objective function. The removed operations are stored in a list of unscheduled operations, following the order in which the destroy operator extracted them from $\mathtt{s}$. We consider two destruction operators, described in the following:

\begin{enumerate}
    \item \texttt{RandomDestroy}: Extracts operations randomly, as proposed by \citet{ruiz2007simple}. 
    
    \item \texttt{PseudoRandomDestroy}: Splits operations into two groups. The first group containing operations with a higher impact in a given solution and the second containing the remaining ones. The impact of each operation $O_i \in \cO$ is computed as $f(\mathtt{s})-f'(\mathtt{s})$, where $f(\mathtt{s})$ and $f'(\mathtt{s})$ are the objective function values of solution~$\mathtt{s}$ before and after removing operation~$O_i$. The operations' impacts are updated iteratively each time an operation is removed. In this approach $d/2$ operations are randomly selected from each group. This operator is inspired by the work of \citet{ruiz2018iteratedgreedy}.
    
\end{enumerate}

We also consider two operators for reconstructing solutions, executed within the \texttt{Repair} procedure~(Line 8 of Algorithm~\ref{alg:ils}). Both operators pick the unscheduled operations one-by-one according to the list order defined during the destroying phase to reinsert in the solution until a new complete solution is created. The repair operators are described as follows:

\begin{enumerate}
    \item \texttt{GreedyRepair}: Uses a greedy rule, computing the objective function value according to the picked operation's insertion in all possible positions between all eligible machines, assigning it to the best possible position.

    \item \texttt{PseudoGreedyRepair}: Reinserts the first $d/2$ operations within the unscheduled list following the first operator's greedy rule. Then, the remaining operations are reinserted at random in any position within the subset of eligible machines for each operation. This operator is inspired by the work of \citet{ruiz2018iteratedgreedy}. 

\end{enumerate}

\begin{figure}[htb] 
	\centerline{\includegraphics[scale=0.18, trim = 0 0 0 0, clip]{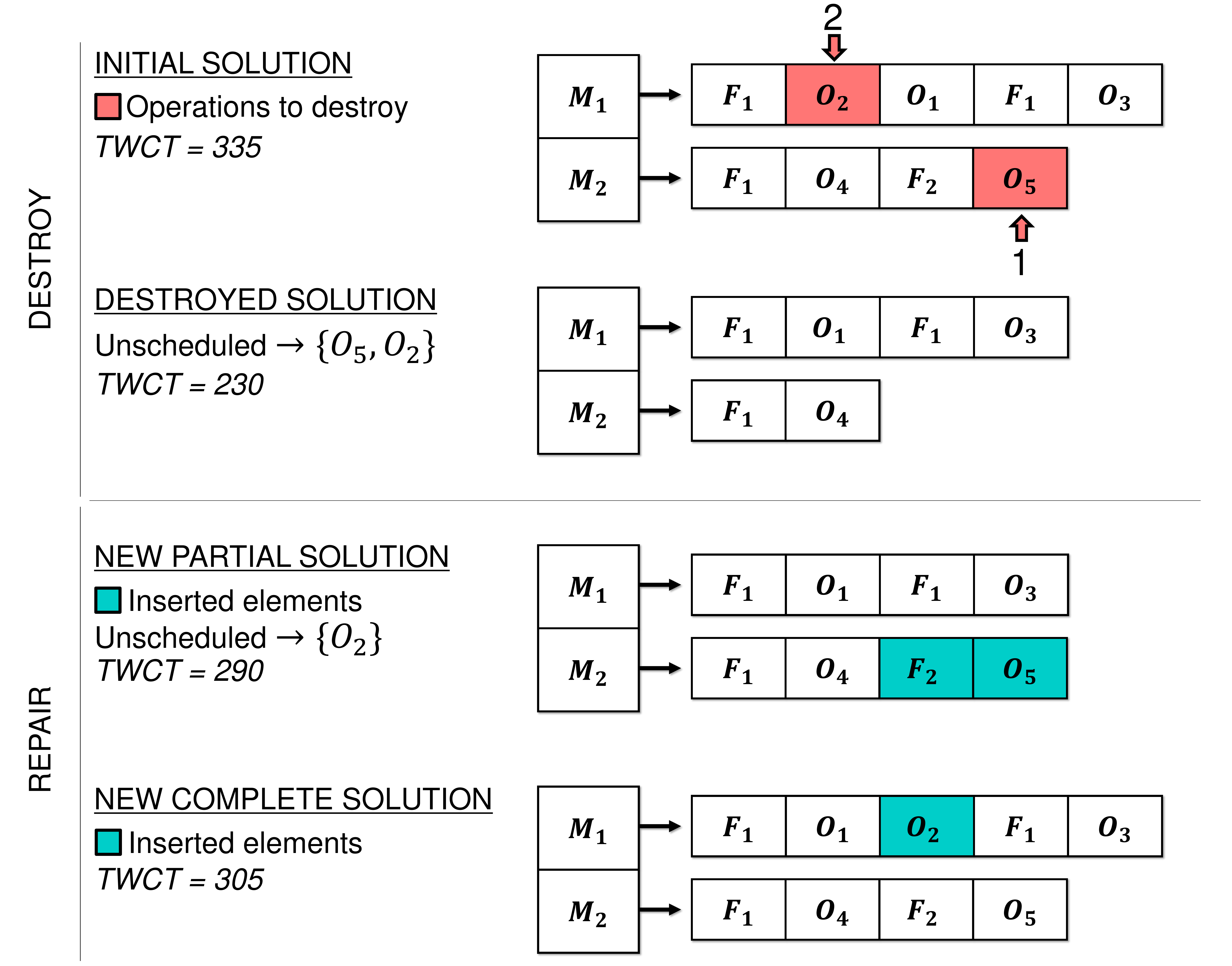}}
	\caption{Destroy and repair exaple with $d=2$. \label{fig:destroy_repair-example}}
\end{figure}

In Figure~\ref{fig:destroy_repair-example}, an example using the \texttt{RandomDestroy} and \texttt{GreedyRepair} operators is shown, considering the example depicted in Figure~\ref{fig:scheduling-example-one}, with $d=2$. Note that operations $O_5$ and $O_2$ are removed from the solution, generating a partial schedule with a TWCT of 230. Then, operation $O_5$ is reinserted in its best position, which, in this case, is its original position. One can note that a setup time is included before operation $O_5$. Finally, operation $O_2$ is inserted in the first batch on machine $M_1$ before operation $O_1$, generating a new complete solution with a TWCT of 305, better than the initial solution.

\subsection{Acceptance Criterion}
\label{sec:accep}

The \texttt{Accept} procedure~(Line 12 of Algorithm~\ref{alg:ils}) receives two solutions, the current solution~$\mathtt{s}$ and a candidate solution $\mathtt{s'}$ to decide whether $\mathtt{s'}$ should replace $\mathtt{s}$ as the current solution. For this, a simulated annealing criterion is employed, accepting $\mathtt{s'}$ with probability $e^{- \Delta / \tau}$, where $\Delta=f(\mathtt{s'})-f(\mathtt{s})$, and $\tau$ is the current temperature. The temperature~$\tau$ starts with a value $\tau_0$, decreasing at each iteration as $\tau = \tau  \kappa$, where $\kappa=[0,1)$ is the cooling rate~\citep{kirkpatrick-1983}. The initial and final temperatures~($\tau_0$ and $\tau_F$) are instance-dependent, as proposed by~\citet{Pisinger2007}, computed as $\tau_0= - (\delta_1  f_o) / \text{ln}(0.5)$ and $\tau_F= - (\delta_2  f_o) / \text{ln}(0.5)$, respectively, where $f_0$ is the cost of the initial solution, and $\delta_1$ and $\delta_2$ are adjustable parameters. The initial and final temperature are set to accept, with 50\% probability, solutions that are $\delta_1$ and $\delta_2$ worse than the initial solution, respectively. Based on these temperatures, the cooling rate~$\kappa$ is set by considering the number of iterations $\eta$ to execute, computed as $\kappa=(\tau_F/\tau_0)^{1/\eta}$. 

\subsection{Infeasibility Strategy}
\label{secmet:deas}

As mentioned earlier, we use an infeasibility strategy to enlarge the problem's search space, increase diversification, and help the algorithm escape from local optimal solutions. In this strategy, the capacity constraint is relaxed, and violating solutions are penalized. Thus, the cost of a solution $\mathtt{s}$ with capacity violations is updated by using a penalty factor $\rho \geq 1$, as $f(\mathtt{s}) = f(\mathtt{s}) + \rho  V$, where $V$ is the total capacity violation among all batches from all machines, and $f(\mathtt{s})$ is the total weighted completion time of solution $\mathtt{s}$. We use two parameters $\rho^{+} \in [0, 1)$ and $\rho^{-} \in [0, 1)$ to update the value of $\rho$ at each iteration. Thus, every time an infeasible solution is accepted, the penalty factor value is increased as $\rho=\rho  (1 + \rho^{+})$. Contrariwise, if the accepted solution is feasible, the factor value is decreased as  $\rho=\rho (1 -\rho^{-})$. The idea is to increase the penalty factor when infeasible solutions are accepted, to force the algorithm to prioritize feasible solutions. The \texttt{Feasible} procedure~(Line 11 of Algorithm~\ref{alg:ils}) updates the parameter~$\rho$ after checking whether a given solution is feasible or not, returning \texttt{true} when no capacity violations exists, and \texttt{false}, otherwise.

\section{Computational Experiments}
\label{sec:experiments}

In this section, we present the conducted computational experiments to evaluate the proposed IG algorithm's performance, organizing them into two parts. First, we present the experiments performed on a set with 72 benchmark instances, developed by \citet{abumarrul2019}. The authors generated these instances based on real data from a ship scheduling problem, with up to 8 machines and 50 operations. In this part, only the most regular destroy and repair operators are considered (\texttt{RandomDestroy} and \texttt{GreedyRepair}). The algorithm's input parameters are calibrated, and the method is compared with the state-of-the-art algorithms for solving the studied problem. In the second part, four IG variants are tested within a new set of instances that combine concepts from different machine scheduling works in the literature with up to 100 operations and 10 machines. We call the instance sets \textit{Small-Sized Benchmark Set} and \textit{Large-Sized Benchmark Set}, respectively. We use C\texttt{++} language for coding the IG algorithm, and ten independent runs are performed for each variant within the experiments. All experiments are performed on a computer with an Intel i7-8700K CPU of 3.70GHz and 64 GB of RAM, running Linux with a single thread. 

In all analyses, we evaluate solutions in terms of the Relative Percentage Deviation~(RPD) concerning the best solutions found for each instance, computed according to Equation~\eqref{eq:RPD}. $TWCT^{Sol}$ denotes the total weighted completion for a given solution of a specific instance, and $TWCT^{Best}$ designates the total weighted completion time regarding the best solution found for the same instance in a given experiment.

\begin{equation}\label{eq:RPD}%
  \mathit{RPD} =%
  \frac{TWCT^{Sol} - TWCT^{Best}}{TWCT^{Best}} \times 100.
\end{equation}

\subsection{Experiments on the Small-Sized Benchmark Set}
\label{sec:small_bench}

\subsubsection{Instances Description}
\label{sec: instances}

As mentioned before, the first benchmark set is composed of 72 instances, proposed by~\citet{AbuMArrul2020}, and available online~\citep{abumarrul2019}. The number of  machines is defined as $m = \{4, 8\}$, and the number of operations as $o=\{15, 25, 50\}$. The number of jobs is defined as $n=\lfloor o/3\rfloor$ and the number of families as $f=3$. The authors consider three input parameters to generate instances with different ranges for the release dates, eligibility levels, and probabilities of associating jobs to operations. The operation's processing times were drawn from a discrete uniform distribution $U(1,30)$ and the operation's sizes from another discrete uniform distribution $U(0,100)$ with step size 10. Note that, in this set, processing times are low. Moreover, the number of operations within a batch may be small, given the possibility of operations with larger sizes close to the machine's capacities, ranging from 80 to 100. We refer the reader to the work of \citet{abumarrul2019} for a complete description of the instances generation process.

\subsubsection{Results and Calibration}
\label{sec:results_small}

Before comparing the IG algorithm with other methods for solving the problem, we set the algorithm's input parameters following a two-phase tuning strategy introduced by \citet{ropke2006adaptive}. The first is a trial-and-error phase, in which the parameters are defined during the algorithm development. The second is the improvement phase, in which a fine-tuning is performed in each parameter individually, within pre-defined possible values, while the remaining parameter values are fixed. We considered the following parameters order for the improvement phase: (1)~Initial temperature parameter for the simulated annealing~($\delta_1$), with values ranging from 0.3 to 0.7 and a step size of 0.1; (2)~Final temperature parameter for the simulated annealing~($\delta_2$) with the possible values of $10^{-3}$, $10^{-4}$, $10^{-5}$, and $10^{-6}$. (3)~Perturbation parameter~($\varepsilon$) with values ranging from 0.05 to 0.20 and a step size of 0.05; (4)~Solution restore parameter~($\lambda$), within the same range of values tested for $\varepsilon$. (5)~Penalty update factor when infeasible solutions are reached~($\rho^{+}$), defined within the range 0.10 to 0.25, with a step size of 0.05; (6)~Penalty update factor when feasible solutions are reached~($\rho^{-}$) defined within the range of 0.05 to 0.20, with a step size of 0.05. All parameters were tuned using the two most regular destroy and repair operators considered in the literature for IG algorithms, the \texttt{RandomDestroy} operator and the \texttt{GreedyRepair} operator. We named \texttt{IG-RG} this variant of the method. Ten independent runs on each instance were performed with $\eta$ = 2500 iterations, which is the value defined in the preliminary tests during the trial-and-error phase. The final parameter values, shown in Table~\ref{tab:parametersdef}, are the ones that presented the best average relative percentage deviation ($\overline{RPD}$) values concerning the best solutions achieved at each step of the parameterization.

\begin{table}[htbp]
	\centering
	\footnotesize
	\caption{Final parameter values for the proposed algorithm.}	
	\begin{tabularx}{\textwidth}{cXcc}
	\hline
		Parameter & Description & Domain & Value \\
         \hline
         $\delta_1$ & Initial temperature definition parameter for the simulated annealing criterion & [0, 1] & 0.6 \\
         $\delta_2$ & Final temperature definition parameter for the simulated annealing criterion & [0, 1] & $10^{-5}$ \\
         $\varepsilon$ & Proportion of the total number of operations to destroy at each iteration & [0, 1] & 0.15 \\
         $\lambda$ & Restore solution parameter & [0, 1] & 0.1 \\
        $\rho^{+}$ & Parameter to update the penalty factor when infeasible solutions are accepted & [0, 1) & 0.20 \\
        $\rho^{-}$ & Parameter to update the penalty factor when feasible solutions are accepted & [0, 1) & 0.05 \\
		\hline			
	\end{tabularx}%
	\label{tab:parametersdef}%
\end{table}%

In Table \ref{tab:components}, we evaluate the impact in the solution quality, in terms of the average relative percentage deviation ($\overline{RPD}$), standard deviation of the RPDs (SD), and the average computational time ($\overline{Time}$), in seconds, with the removal of each feature of the algorithm individually. Each configuration corresponds to the disabling of one of the following features: (LS)~RVND local search; (SA)~Simulated annealing acceptance criterion; (DR)~Destroy and Repair steps; (Inf.)~Infeasibility strategy; (Rest.)~Solution restore strategy. When the destroy and repair steps are turned off, we replace it with a regular perturbation step in which we select one neighborhood at random, and $d = \lceil \varepsilon o \rceil$ random moves are performed. Note that the algorithm using the perturbation instead of the destruction and repair operators follows an Iterated Local Search (ILS) metaheuristic structure. As one can notice, the local search is the most relevant component of the method and the most time-consuming feature. The algorithm's average computational time reduces when it disregards the diversification strategies (simulated annealing and infeasibility strategy). However, the small addition of time is worth it due to the deviation reduction these strategies bring to the method. We can see that the configuration with all features combined generates the smallest $\overline{RPD}$ and SD, indicating that they all contribute to the algorithm's performance. 

\begin{table}[htbp]
  \centering
  \caption{Average RPD and computational time with the removal of each algorithm's feature.}
    \begin{tabular}{lrcccccccccc}
    \toprule
    \multicolumn{1}{c}{Config.} &       & LS    & SA    & DR    & Inf.  & Rest. &       & $\overline{RPD}$ & SD &       & $\overline{Time}$ \\
    \midrule
    No LS &       &       & \tiny $\bullet$      & \tiny $\bullet$      & \tiny $\bullet$      & \tiny $\bullet$      &       & 2.68  & 2.38  &       & 0.15 \\
    No SA &       & \tiny $\bullet$      &       & \tiny $\bullet$      & \tiny $\bullet$      & \tiny $\bullet$      &       & 0.22  & 0.42  &       & 7.70 \\
    No DR &       & \tiny $\bullet$      & \tiny $\bullet$      &       & \tiny $\bullet$      & \tiny $\bullet$      &       & 0.37  & 0.90  &       & 12.86 \\
    No Inf. &       & \tiny $\bullet$      & \tiny $\bullet$      & \tiny $\bullet$      &       & \tiny $\bullet$      &       & 0.19  & 0.31  &       & 7.69 \\
    No Rest. &       & \tiny $\bullet$      & \tiny $\bullet$      & \tiny $\bullet$      & \tiny $\bullet$      &       &       & 0.18  & 0.32  &       & 8.20 \\
    Complete &       & \tiny $\bullet$      & \tiny $\bullet$      & \tiny $\bullet$      & \tiny $\bullet$      & \tiny $\bullet$      &       & 0.15  & 0.26  &       & 8.18 \\
    \bottomrule
    \end{tabular}%
  \label{tab:components}%
\end{table}%

\begin{figure}[htbp!] 
	\centerline{\includegraphics[scale=0.65, trim = 0 0 0 0, clip]{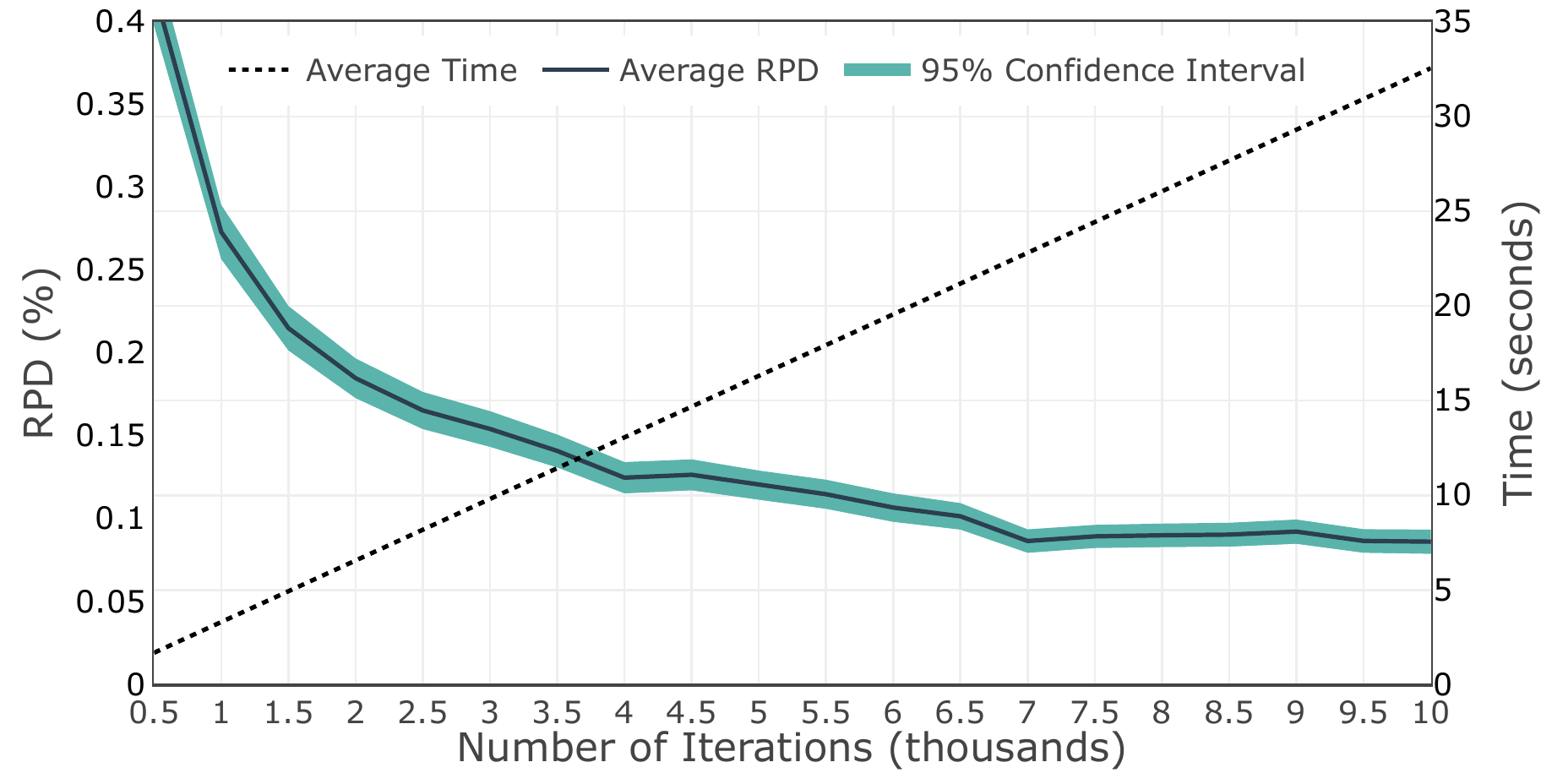}}
	\caption{Average RPD and average computational time with different number of iterations. \label{fig:scheduling-example}}
\end{figure}

In Figure \ref{fig:scheduling-example}, we evaluate the trade-off between the solution quality in terms of the $\overline{RPD}$ and the average computational time, regarding the number of iterations considered. We run the experiments with the total number of iterations~($\eta$) ranging from 500 to 10000 with a step size of 500. Note that the average time grows linearly as the number of iterations increases. The $\overline{RPD}$ is presented with a 95\% confidence interval. One can note that the average RPD decreases as the number of iterations grows, reducing faster up to 4000 iterations. From 4000 to 4500, the method is stable. It returns to a continuous reduction of the average RPD, but slower, from 4500 to 7000 iterations. Then, from 7000 iterations, the method stabilizes, maintaining the same quality of the solutions until reaching the maximum number of iterations tested~(10000), indicating that it converges to a $\overline{RPD}$ around 0.08\%. 

In the next analysis, we compare the \texttt{IG-RG} with the best methods in the literature for solving the problem. As mentioned before, in \citet{abu2020matheuristics}, the authors tackled the problem with matheuristic algorithms, comparing six variants of it. In their analysis, two of the matheuristics presented the best results, the \texttt{GRASP-Math}$_3$ and the \texttt{ILS-Math}$_3$. Both methods use a combination of a metaheuristic approach with MIP-based local searches using a batch scheduling formulation. The former considers a Greedy Randomized Adaptive Search Procedure~(GRASP) framework, and the latter an Iterated Local Search~(ILS). We compare them in terms of the RPD, computational time, and capacity of reaching the best-found solution. Based on the previous analysis, we ran the algorithm for 2500, 4500, and 7000 iterations. The results are depicted in Table~\ref{tab:bench1_math_comparison} with the average~(Avg), maximum~(Max) and standard deviation~(SD) for the RPDs and computational times, the percentage of instances~(\%Inst) and runs~(\%Runs) in which each algorithm achieved the best solution, and the percentage of instances in which the best solution is uniquely found by a given approach (\%Unique). Each \texttt{IG-RG} variant includes the number of iterations in its name. Our experiments were performed on a machine with the same characteristics as the one used by \citet{abu2020matheuristics} for a fair comparison.  Note that the \texttt{IG-RG} has a superior performance in all criteria also when running for 2500 iterations (\texttt{IG-RG-2500}). Its worst-case RPD (Max) remains around 2\%, with a low standard deviation (0.28), showing consistency within the different runs. Even when the number of iterations grows, as in the case of the \texttt{IG-RG-7000} (\texttt{IG-RG} running for 7000 iterations), the average computational time is at least 94\% lower than the matheuristic ones. The \texttt{IG-RG} provided new best solutions (upper bound) for this set of instances. The algorithm reached the best solution in 94.44\% of the instances and 70\% of the runs when executed by 7000 iterations. Complete detailed results are available as supplementary material.

\begin{table}[htbp]
  \centering
  \caption{Comparison between the iterated greedy algorithm and the best methods in the literature for this set of instances.}
  \resizebox{\textwidth}{!}{%
    \begin{tabular}{lrccccccccccc}
    \toprule
    \multicolumn{1}{c}{\multirow{2}[4]{*}{Algorithm}} &       & \multicolumn{3}{c}{RPD (\%)} &       & \multicolumn{3}{c}{Time (seconds)} &       & \multicolumn{3}{c}{Achieved Best Solution} \\
\cmidrule{3-5}\cmidrule{7-9}\cmidrule{11-13}          &       & Avg   & Max   & SD    &       & Avg   & Max   & SD    &       & \%Inst. & \%Run & \%Unique \\
    \midrule
    \texttt{GRASPMath}$_3$ &       & 0.79  & 20.57 & 1.28  &       & 464.68 & 2677.43 & 529.83 &       & 62.50 & 40.56 & 0.00 \\
    \texttt{ILS-Math}$_3$ &       & 0.74  & 6.62  & 0.96  &       & 411.67 & 2735.17 & 503.52 &       & 58.33 & 37.64 & 1.39 \\
    \texttt{IG-RG-2500} &       & 0.15  & 2.03  & 0.28  &       & 8.18  & 26.78 & 9.71  &       & 73.61 & 61.25 & 0.00 \\
    \texttt{IG-RG-4500} &       & 0.11  & 1.47  & 0.23  &       & 14.69 & 47.28 & 17.48 &       & 83.33 & 64.86 & 4.17 \\
    \texttt{IG-RG-7000} &       & 0.07  & 1.23  & 0.16  &       & 22.80 & 75.26 & 27.14 &       & 94.44 & 70.00 & 19.44 \\
    \bottomrule
    \end{tabular}%
    }
  \label{tab:bench1_math_comparison}%
\end{table}%

\subsection{Experiments on the Large-Sized Benchmark Set}
\label{sec:large_bench}

\subsubsection{Instances Description}
\label{sec:large_instances}

The second benchmark combines concepts from different instance generation processes in the scheduling literature. We generated 3 instances for each combination of the number of machines $m = \{5, 10\}$, the number of operations $o=\{50, 75, 100\}$, the number of jobs $n=o/q$, where $q=\{3,5\}$, and the number of families $f=\{3,5\}$. Therefore, the set is composed of 72 instances in total. The generation of the remaining aspects is described in the following.

Processing times follow a discrete uniform distribution $U(1,100)$. This is a well-know interval considered in the scheduling literature for drawing processing times (see, for instance, the works of \citet{sheikhalishahi2019multi}, \cite{fanjul2010iterated}, \cite{lin2013multiple}, and \cite{akturk2001new}). Family setup times follow a discrete uniform distribution $U(10,30)$, as in the work of \citet{sheikhalishahi2019multi}. And, job's weights are drawn from a discrete uniform distribution $U(1,10)$, as established by \citet{lin2013multiple}. 

Operation's families are defined according to the scheme proposed by \citet{dunstall2005comparison}. First, a random number between 0 and 1 is generated for each family. Then, these numbers are weighted to define the proportion of operations in each family. Following the defined proportion, operations are assigned to families at random.

Eligibility subsets are defined according to the scheme proposed by \citet{bitar2016memetic}. First, a number of eligible machines for each operation is drawn from a discrete uniform distribution $U(1,m)$. Then, the eligible machines are selected at random, respecting the defined number.

Release dates for machines and operations are defined based on the works of \citet{pei2020new}, \cite{AbuMArrul2020}, and \cite{akturk2001new}, following a discrete uniform distribution $U(0,MR)$, where $MR$ is the maximum release date, defined as $MR = \bigl\lceil \alpha  \bigl( \sum_{O_i \in \cO}{p_i + \sum_{g \in \cF} s_{g}  |\cF_g|\bigr)}/ m \bigr\rceil$. We use $\alpha=0.5$ for defining the operation's release dates and $\alpha=0.1$ for defining the machine's release dates.

The association between jobs and operations follows a modified version of the scheme proposed by \citet{AbuMArrul2020}, and we organized it into two steps. Let $\cU \subseteq \cO$ be a subset of unassigned operations, initialized as $\cU = \cO$. In the first step, jobs are selected one by one at random, and $\max\{|\cU|, \lceil o/n \rceil\}$ operations are randomly associated with it, being selected from subset $\cU$. Selected operations are removed from $|\cU|$ before moving to the next job. After the first step, each operation is associated with a unique job. Then, in the second step, for each pair operation/job, we define if they are associated with 5\% of probability. 

Capacities of machines and operation's sizes are defined based on the works of \citet{pei2015serial} and \citet{AbuMArrul2020}. The machine's capacities follow a discrete uniform distribution $U(10,15)$, while the operation's sizes are drawn from a discrete uniform distribution $U(1,5)$.

\subsubsection{Results and Discussion}
\label{sec: results}

After highlighting the superior performance of the proposed IG algorithm against the current methods in the literature, in this section, we conduct experiments within the new set of instances considering four variants of the method, with different combinations among the destroy and repair operators (Section~\ref{sec:destroy_repair}): (1)~\texttt{IG-RG}: IG with \texttt{RandomDestroy} and \texttt{GreedyRepair}; (2)~\texttt{IG-RP}: IG with \texttt{RandomDestroy} and \texttt{PseudoGreedyRepair}; (3)~\texttt{IG-PG}: IG with \texttt{PseudoRandomDestroy} and \texttt{GreedyRepair}; (4)~\texttt{IG-PP}: IG with \texttt{PseudoRandomDestroy} and \texttt{PseudoGreedyRepair}. The results are shown in Table~\ref{tab:bench2_methods}, with instances grouped by the number of operations~($o$) and machines~($m$). In this experiment, we diversify the total number of iterations~($\eta$) to evaluate the quality of the solutions in terms of the average relative percentage deviation~($\overline{RPD}$) and the average computational time~($\overline{Time}$) for 2500, 4500 and 7000 iterations. Best values for the $\overline{RPD}$ are highlighted in bold. 

\begin{table}[htb!]
  \centering
  \caption{Results by group of instances for each algorithm within different number of iterations.}
  \resizebox{\textwidth}{!}{%
    \begin{tabular}{cccccccccccccc}
    \toprule
    \multirow{2}[4]{*}{$\eta$} & Group &       & \multicolumn{2}{c}{\texttt{IG-RG}} &       & \multicolumn{2}{c}{\texttt{IG-RP}} &       & \multicolumn{2}{c}{\texttt{IG-PG}} &       & \multicolumn{2}{c}{\texttt{IG-PP}} \\
\cmidrule{4-5}\cmidrule{7-8}\cmidrule{10-11}\cmidrule{13-14}          & $o$-$m$ &       & $\overline{RPD}$   & $\overline{Time}$  &       & $\overline{RPD}$   & $\overline{Time}$  &       & $\overline{RPD}$   & $\overline{Time}$  &       & $\overline{RPD}$   & $\overline{Time}$ \\
    \midrule
    \multirow{7}[4]{*}{2500} & 50-5  &       & 0.35  & 16.79 &       & 0.52  & 38.61 &       & \textbf{0.34} & 17.37 &       & 0.54  & 40.14 \\
          & 50-10 &       & \textbf{0.46} & 17.62 &       & 0.60  & 31.93 &       & 0.48  & 18.19 &       & 0.61  & 33.19 \\
          & 75-5  &       & \textbf{0.69} & 74.81 &       & 1.32  & 182.71 &       & 0.74  & 75.80 &       & 1.41  & 195.91 \\
          & 75-10 &       & \textbf{0.73} & 68.35 &       & 1.46  & 144.41 &       & 0.76  & 70.07 &       & 1.56  & 150.12 \\
          & 100-5 &       & \textbf{0.72} & 216.27 &       & 1.70  & 583.66 &       & 0.81  & 216.28 &       & 1.88  & 629.27 \\
          & 100-10 &       & \textbf{0.75} & 189.12 &       & 1.85  & 441.93 &       & 0.84  & 194.83 &       & 1.98  & 468.39 \\
\cmidrule{2-14}          & All   &       & \textbf{0.62} & 97.16 &       & 1.24  & 237.21 &       & 0.66  & 98.76 &       & 1.33  & 252.83 \\
    \midrule
    \multirow{7}[4]{*}{4500} & 50-5  &       & \textbf{0.26} & 30.15 &       & 0.38  & 69.47 &       & 0.28  & 31.28 &       & 0.41  & 72.43 \\
          & 50-10 &       & \textbf{0.38} & 31.69 &       & 0.44  & 57.32 &       & 0.41  & 32.56 &       & 0.50  & 59.77 \\
          & 75-5  &       & \textbf{0.49} & 134.17 &       & 1.20  & 329.93 &       & 0.59  & 134.84 &       & 1.31  & 351.52 \\
          & 75-10 &       & \textbf{0.55} & 122.83 &       & 1.33  & 259.57 &       & 0.56  & 126.27 &       & 1.45  & 270.02 \\
          & 100-5 &       & \textbf{0.52} & 388.95 &       & 1.57  & 1050.79 &       & 0.60  & 390.18 &       & 1.69  & 1138.61 \\
          & 100-10 &       & \textbf{0.61} & 341.16 &       & 1.72  & 796.96 &       & \textbf{0.61} & 350.29 &       & 1.86  & 841.75 \\
\cmidrule{2-14}          & All   &       & \textbf{0.47} & 174.82 &       & 1.11  & 427.34 &       & 0.51  & 177.57 &       & 1.20  & 455.69 \\
    \midrule
    \multirow{7}[4]{*}{7000} & 50-5  &       & \textbf{0.19} & 47.33 &       & 0.31  & 109.10 &       & 0.23  & 49.20 &       & 0.33  & 113.87 \\
          & 50-10 &       & \textbf{0.30} & 49.99 &       & 0.39  & 90.57 &       & 0.33  & 51.30 &       & 0.40  & 94.11 \\
          & 75-5  &       & \textbf{0.45} & 209.51 &       & 1.10  & 515.05 &       & 0.47  & 211.97 &       & 1.21  & 550.21 \\
          & 75-10 &       & 0.48  & 192.35 &       & 1.24  & 406.24 &       & \textbf{0.47} & 199.05 &       & 1.31  & 423.76 \\
          & 100-5 &       & \textbf{0.42} & 610.52 &       & 1.49  & 1640.53 &       & 0.45  & 608.87 &       & 1.64  & 1778.70 \\
          & 100-10 &       & \textbf{0.48} & 536.75 &       & 1.62  & 1249.20 &       & 0.50  & 550.54 &       & 1.75  & 1322.61 \\
\cmidrule{2-14}          & All   &       & \textbf{0.39} & 274.41 &       & 1.02  & 668.45 &       & 0.41  & 278.49 &       & 1.11  & 713.88 \\
    \bottomrule
    \end{tabular}%
  }
  \label{tab:bench2_methods}%
\end{table}%

One can see that methods with the \texttt{GreedyRepair} operator (\texttt{IG-RG} and \texttt{IG-PG}) outperformed the ones with the \texttt{PseudoGreedyRepair} operator (\texttt{IG-RP} and \texttt{IG-PP}), in terms of the $\overline{RPD}$ and $\overline{Time}$, independently of the number of iterations executed. Moreover, the discrepancy between these algorithms grows when more iterations are performed. For example, when running for 7000 iterations and considering the complete set of instances, \texttt{IG-RG} reached an average RPD of 0.39\%. At the same time, algorithms with the \texttt{PseudoGreedyRepair} operator remained with average RPDs above 1\%. Not surprisingly, the $\overline{RPD}$ values decrease when the number of iterations increases for all algorithm variants. One can also note that the algorithms with the \texttt{GreedyRepair} operator consume less than half the others' computational time. Overall, the \texttt{IG-RG} dominates the remaining variants with the best values for the $\overline{RPD}$ and $\overline{Time}$ on almost all groups, independently of the number of iterations performed. An interesting aspect of the results is how the average computational time is usually higher for groups with five machines. This is mainly because solving these instances is more challenging due to the higher number of operations performed by each machine, increasing the local search's impact on the algorithm's performance. The studied problem is characterized by the fact that few operations define the job's completion times and, consequently, the objective function's value. Thus, any random insertion of operations can cause a significant disruption in a given solution's cost. The local search may struggle to restore reasonable solutions from a perturbed solution, given that numerous movements lead to solutions of equivalent cost. This corroborates the weak performance of the \texttt{PseudoGreedyRepair} operator in the addressed problem.

\begin{figure}[h!]
\centering
    \includegraphics[scale=0.55, trim = 0 0 0 0,clip]{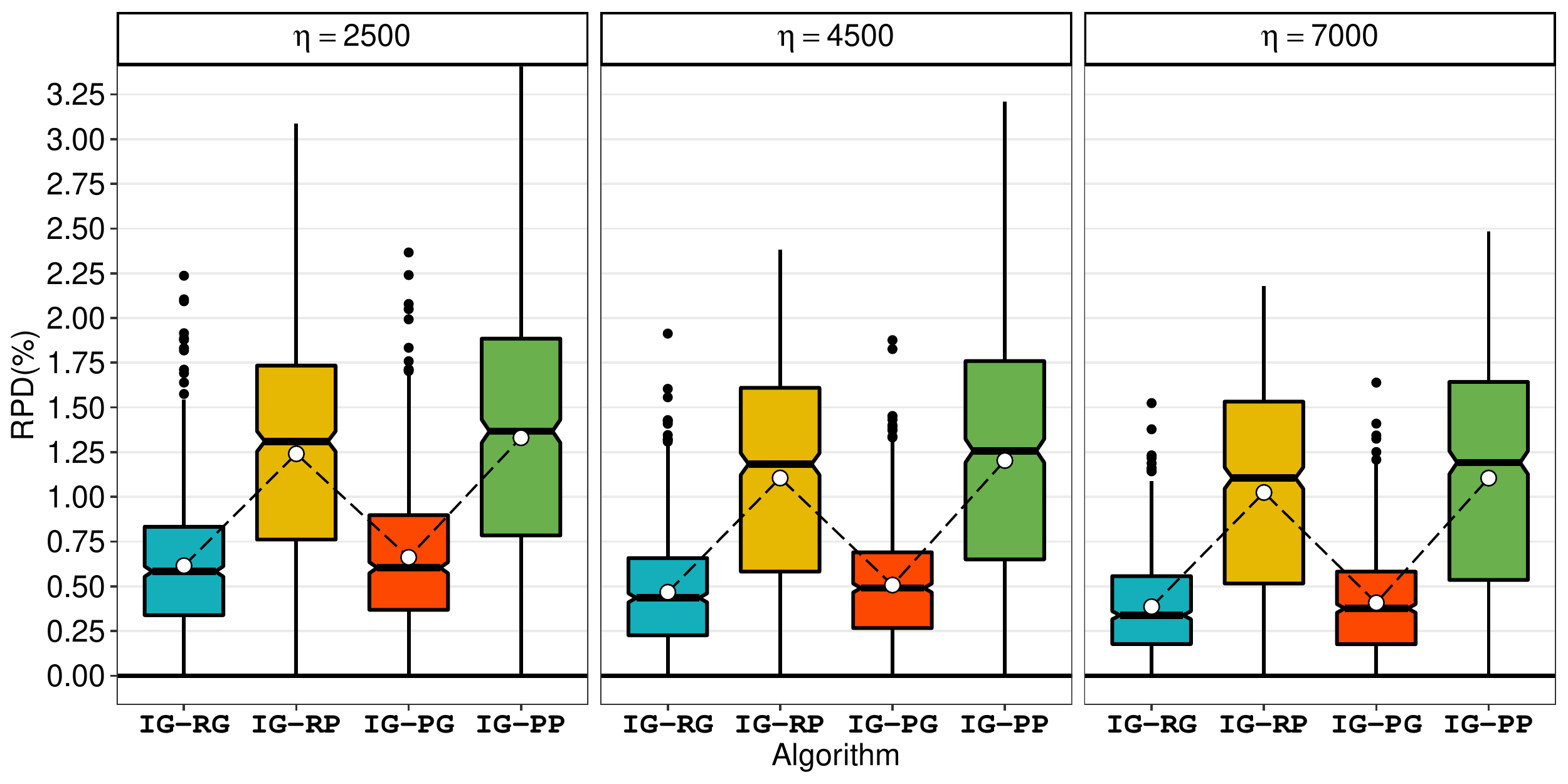}
    \caption{Boxplots of the RPD distributions for each algorithm within different number of iterations, considering the complete set of instances.}
    \label{fig:boxplot_operations}
\end{figure}

Figure~\ref{fig:boxplot_operations} shows the RPD distributions of each algorithm for the different numbers of iterations considering the complete set of instances to help visualize the contrast between the methods' solutions. One can note some relevant points: (1)~The distributions get less spread out and with smaller interquartile values as the number of iterations grow. (2)~The weakness of the algorithms that use the \texttt{PseudoGreedyRepair} operator is more evident, getting worse when the operator is combined with the \texttt{PseudoRandomDestroy} operator. (3)~Faster improvement by the \texttt{IG-RG} algorithm can be noted, with it reaching an average RPD below 0.5\% when running for 4500 iterations. (4)~Slower convergence of the \texttt{IG-PG} variant compared to \texttt{IG-RG} but competitive when performed for 7000 iterations.

To improve the discussion and provide more statistical information over the RPD distributions, we conducted the pairwise Wilcoxon rank-sum test with Hommel's \textit{p}-values adjustment with a confidence level of 0.05. Before performing the test, we executed the Shapiro-Wilk test on the distributions, which indicated that the RPDs are not normally distributed. The \textit{p}-values obtained from the pairwise Wilcoxon test were all below 0.05, except when comparing the \texttt{IG-RG} against \texttt{IG-PG}, running for 7000 iterations, where the \textit{p}-value is 0.16. This result confirms what we highlighted, and shows that significant statistical differences can be noted between the algorithms using the \texttt{GreedyRepair} (\texttt{IG-RG} and \texttt{IG-PG}) and \texttt{PseudoGreedyRepair} operators (\texttt{IG-RP} and \texttt{IG-PP}). Also, although \texttt{IG-RG} achieves a better value for $\overline{RPD}$ regardless of the number of iterations performed, the statistic test does not indicate a significant difference between it and the \texttt{IG-PG} variant when the algorithms are executed for 7000 iterations, showing that both methods converge for solutions of similar quality as the number of iterations grows. 

In the next analysis, we evaluate the methods' performance in terms of the proportion of best solutions achieved by each one. Thus, Table~\ref{tab:percentage_instances_bench2} depicts the percentage of instances (\%Inst.) and runs (\%Run) in which each method achieved the best solution and the percentage of instances in which the best solution was uniquely achieved by a given algorithm (\%Unique). The table includes the average RPD ($\overline{RPD}$) and the standard deviation (SD) of the RPDs regarding the runs in which the best solution was not achieved. Algorithm names include the number of iterations executed, and the best results are highlighted in bold. Note that variant \texttt{IG-PG} surpasses the others, regardless of the number of iterations, finding unique best solutions in 38.89\% of the instances (29,17\% + 6.94\% + 2.78\%). Moreover, it reaches the best solutions in 44.44\% of the instances when running for 7000 iterations. These results indicate that the \texttt{PseudoRandomDestroy} operator can lead to better solutions when combined with the \texttt{GreedyRepair} operator. However, variant \texttt{IG-PG} needs a higher number of iterations to achieve an overall competitive performance, which may indicate that the method is more likely to get stuck in local optimum solutions. However, \texttt{IG-PG} seems to be a good option when the decision-maker has more time to solve the problem. One can note that variant \texttt{IG-PP} had the worst performance in these criteria, not providing any of the best solutions. \texttt{IG-RP} variant found some of the best solutions, collaborating with future studies related to this set of instances. Another important aspect concerns the low number of runs in which the methods found the best solution. Note that \texttt{IG-RG-7000} found the best solution in 7.92\% of the runs. In contrast, for the smallest set of instances, this number was 70.00\% (see Section \ref{sec:results_small}), proving the new benchmark to be a more challenging set to be solved. Despite this behavior, the \textit{Not Achieved} columns show that $\overline{RPD}$ and SD are low when the method does not reach the best solutions, highlighting its stable performance.

\begin{table}[htb!]
  \centering
  \caption{Analysis of the best solutions achievement among the IG algorithm variants, considering the complete set of instances.}
    \begin{tabular}{lrrrrrrr}
    \toprule
    \multicolumn{1}{c}{\multirow{2}[4]{*}{Algorithm}} &       & \multicolumn{1}{c}{\multirow{2}[4]{*}{\%Inst.}} & \multicolumn{1}{c}{\multirow{2}[4]{*}{\%Run}} & \multicolumn{1}{c}{\multirow{2}[4]{*}{\%Unique}} &       & \multicolumn{2}{c}{Not Achieved} \\
\cmidrule{7-8}          &       &       &       &       &       & \multicolumn{1}{c}{$\overline{RPD}$} & \multicolumn{1}{c}{SD} \\
    \midrule
    \texttt{IG-RG-2500} &       & 13.89 & 2.64  & 5.56  &       & 0.63  & 0.36 \\
    \texttt{IG-RG-4500} &       & 20.83 & 4.58  & 8.33  &       & 0.49  & 0.30 \\
    \texttt{IG-RG-7000} &       & 40.28 & \textbf{7.92} & 18.06 &       & \textbf{0.42} & \textbf{0.26} \\
    \texttt{IG-RP-2500} &       & 8.33  & 1.11  & 0.00  &       & 1.25  & 0.61 \\
    \texttt{IG-RP-4500} &       & 11.11 & 2.36  & 1.39  &       & 1.13  & 0.60 \\
    \texttt{IG-RP-7000} &       & 13.89 & 3.61  & 1.39  &       & 1.06  & 0.57 \\
    \texttt{IG-PG-2500} &       & 18.06 & 2.50  & 2.78  &       & 0.68  & 0.39 \\
    \texttt{IG-PG-4500} &       & 23.61 & 4.17  & 6.94  &       & 0.53  & 0.31 \\
    \texttt{IG-PG-7000} &       & \textbf{44.44} & \textbf{7.92} & \textbf{29.17} &       & 0.44  & 0.27 \\
    \texttt{IG-PP-2500} &       & 2.78  & 0.42  & 0.00  &       & 1.34  & 0.67 \\
    \texttt{IG-PP-4500} &       & 6.94  & 1.53  & 0.00  &       & 1.22  & 0.65 \\
    \texttt{IG-PP-7000} &       & 11.11 & 2.08  & 0.00  &       & 1.13  & 0.63 \\
    \bottomrule
    \end{tabular}%
  \label{tab:percentage_instances_bench2}%
\end{table}%

\begin{figure}[htbp!] 
	\centerline{\includegraphics[width=0.85\textwidth, trim = 0 40 0 0, clip]{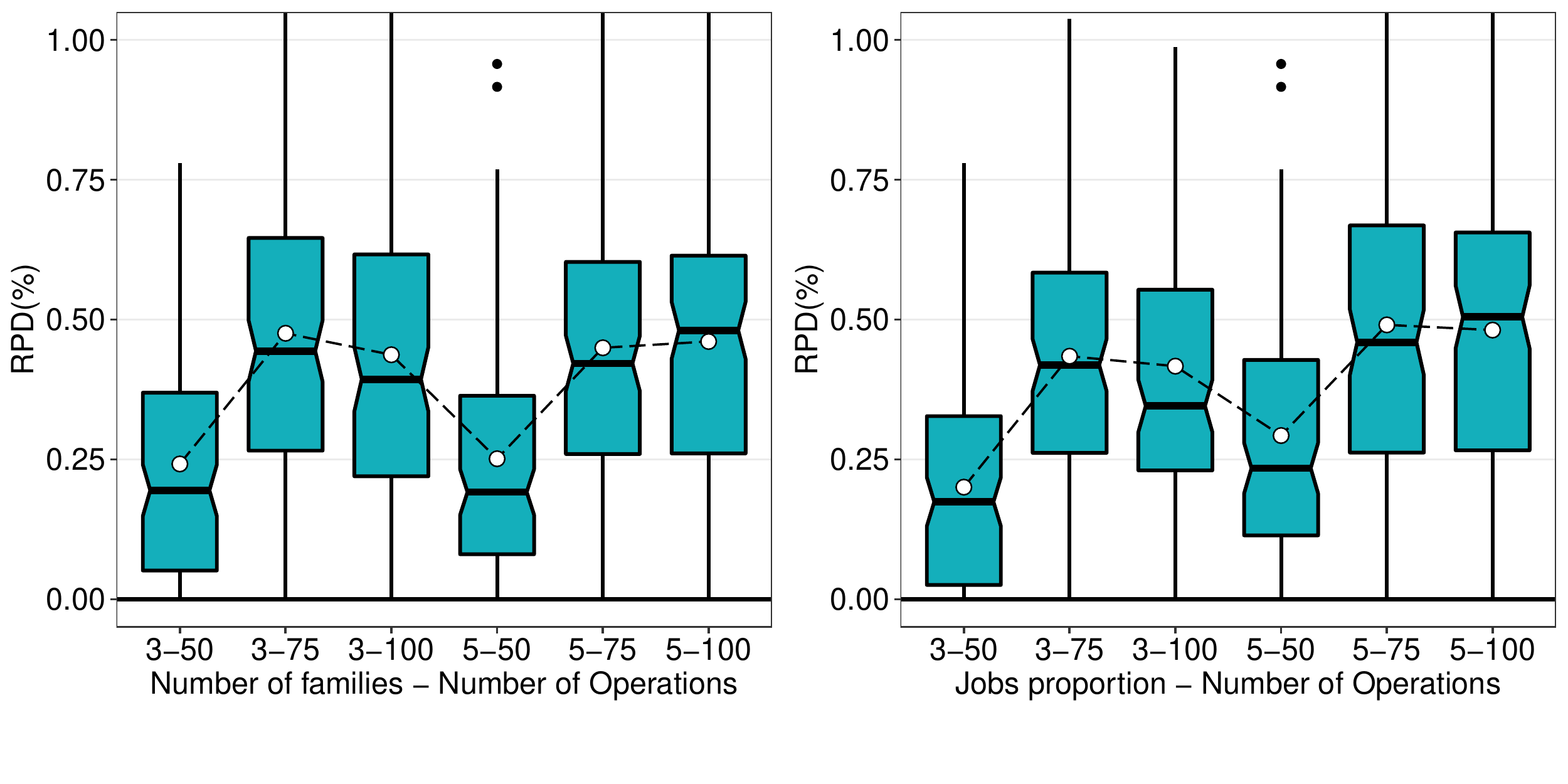}}
	\caption{Boxplot of the RPD distributions for the 
	\texttt{IG-RG} algorithm, running for 7000 iterations, considering different aspects of the instances. \label{fig:specific_parameter_analysis}}
\end{figure}

Figure \ref{fig:specific_parameter_analysis} shows the RPD distributions for two different grouping schemes to assess the performance of the \texttt{IG-RG-7000} concerning relevant aspects of the instances. In the first, we grouped the instances by the number of families~($f$) and operations~($o$). While in the second, the groups indicate the proportion of jobs~($q$), defined in the instance generation process, and the number of operations~($o$) to schedule. One can see an improvement in the solutions when the number of families increases. This behavior indicates that the family constraint reduces the possible combinations of operations within the batches since it concerns a hard constraint, enhancing the local search performance. The same effect can be observed for the number of jobs. When the proportion of jobs is equal to three~($q=3$), meaning that approximately three operations will compose a job, more jobs are considered. For example, when 100 operations must be scheduled, the number of jobs is equal to 34. In these cases, the RPD distributions are less dispersed with lower average and interquartile values. If the number of jobs reduces, fewer operations' completion times define the objective function's value. Thus, the proportion of solutions with equivalent cost increases, limiting the efficiency of the local search.

\section{Conclusions}
\label{sec:conclusion}

In this paper, we addressed a complex parallel machine scheduling problem, where jobs are composed of operations that are arranged into families. Non-anticipatory setup times incur at the beginning of each batch, formed by a sequence of operations considering their sizes and families. The problem also considers release dates for operations and machines, and machine's eligibility and capacity constraints. To solve it, we developed four Iterated Greedy (IG) algorithm variants, combining two destroy operators and two repair operators with a RVND local search procedure and other diversification and intensification strategies.

We tested the algorithms in two instance sets, named \textit{Small-Sized Benchmark Set} and \textit{Large-Sized Benchmark Set}. The former is a known benchmark set from the literature, introduced by \citet{AbuMArrul2020}. The latter is a new set proposed in this work, based on instance generation definitions of several papers from the machine scheduling literature. The results showed that two variants using a greedy repair operator performed better than the remaining ones. The best variant outperformed the current literature approaches for solving the addressed problem in terms of the average deviation from the best solutions and the average computational time. Moreover, the IG algorithm provided new best solutions (upper bounds) in the \textit{Small-Sized Benchmark Set}. Experiments showed that the solution quality improves when the number of iterations of the algorithm grows. Nevertheless, the algorithms are still efficient when running for the small number of iterations tested (2500 iterations). Using the most common operators of destruction and repair in IG algorithms' literature, one of the method's variants presented the best results in different analyses presented, emphasizing that these operators are still relevant and efficient despite their simplicity. 

The present study reinforces the applicability of iterated greedy algorithms in solving combinatorial optimization problems even when a complicated problem inspired by a real context is considered. The idea of destroying and repairing the solution is a simple yet powerful and efficient technique, as we confirmed in our experiments. Since we are dealing with a complex parallel machine scheduling problem, it is worth emphasizing that the developed IG algorithms can be applied to simplified variants of the problem. This aspect enhances its relevance to the machine scheduling literature, supporting researchers studying similar problems or simplified variants of it. In this work, we consider that all the problem data are deterministic and known in advance. However, considering uncertainties in some aspects of the problem would make it attractive to researchers working with realistic problems. Many works on stochastic scheduling problems consider uncertainty in task processing times. Nevertheless, due to the complexity of the problem addressed, uncertainties can be defined for several aspects. For instance, uncertainty in the size of operations could lead to infeasible solutions concerning the machines' capacities. Moreover, uncertainty in the release dates would impact solutions' costs due to non-anticipatory setup times consideration. Therefore, investigating these aspects and developing tools to deal with the problem's stochastic variants is an exciting theme for future work.

\section*{Acknowledgments}

This research was partially supported by PUC-Rio, by the Coordenação de Aperfeiçoamento de Pessoal de Nível Superior – Brasil~(CAPES) – Finance Code 001, by the Conselho Nacional de Desenvolvimento Cient{\'i}fico e Tecnol{\'o}gico~(CNPq), under grant numbers \mbox{313521/2017-4} and \mbox{315361/2020-4}, and by the Norwegian Agency for International Cooperation and Quality Enhancement in Higher Education~(Diku) – Projct number UTF-2017-four-year/10075.

\bibliographystyle{elsarticle-harv}
%\bibliography{references}
{\linespread{1.3}\selectfont\bibliography{references}}

\end{document}